\documentclass[10pt]{amsart}

\usepackage{amsmath,graphics}
\usepackage{amsfonts,amssymb}

\theoremstyle{plain}
\newtheorem*{theorem*}{Theorem}
\newtheorem*{lemma*} {Lemma}
\newtheorem*{corollary*} {Corollary}
\newtheorem*{proposition*} {Proposition}
\newtheorem{theorem}{Theorem}[section]
\newtheorem{lemma}[theorem]{Lemma}
\newtheorem{corollary}[theorem]{Corollary}
\newtheorem{proposition}[theorem]{Proposition}

\theoremstyle{remark}
\newtheorem*{remark}{Remark}
\newtheorem*{definition}{Definition}

\theoremstyle{definition}

\def \R {\mathbb{R}}

\def \Z {\mathbb{Z}}
\def \C {\mathbb{C}}

\def \Qt{\widetilde{Q}}
\def \Qtl{\widetilde{Q}_\lambda}
\def \Pt{\widetilde{P}}
\def \Ptl{\widetilde{P}_{\lambda}}

\def \tI{\mathcal{I}_n}

\def \zI{\zeta^I}
\def \zJ{\zeta^J}
\def \bn{\begin{enumerate}}
\def \en{\end{enumerate}}
\def \bdm{\begin{displaymath}}
\def \edm{\end{displaymath}}

\def \bp{\begin{proof}}
\def\ep{\end{proof}}
\def\av{\arrowvert}

\def\vd{\vdots}
\def\dd{\ddots}
\def\U+e{(U^+)^e}
\def\Qtd{\dot{\widetilde{Q}}}

\def\trm{\hspace{0.08in}\textrm{for} \hspace{0.08in}}
\def\be{\begin{equation}}
\def\ee{\end{equation}}
\def\mfg{\mathfrak{g}}

\def\mft{\mathfrak{t}}
\def\mfb{\mathfrak{b}^+}
\def\mfbm{\mathfrak{b}^{-}}

\def\mfu{\mathfrak{u}^+}
\def\mfum{\mathfrak{u}^{-}}
\def\Fdot{ F_{\mbox{\boldmath{.}}}}
\def\mcvc{\mathcal{V}_n^C}
\def\mcvb{\mathcal{V}_n^B}
\def\mcvd{\mathcal{V}_{n+1}^D}
\def\hs{\hspace{0.07in}}
\def\mcv{\mathcal{V}}

\def\qhoge{qH^*(OG^e(n))}

\def\qhlg{qH^*(LG(n))}
\def\oge{OG^e(n)}
\def\ogo{OG^o(n)}
\def\lg{LG(n)}
\def\In{\mathcal{I}_n}
\def\Im{\mathcal{I}_{n+1}}
\def\SOo{SO_{2n+1}(\C)}
\def\SOe{SO_{2n+2}(\C)}
\def\Sp{Sp_{2n}(\C)}
\def\dn{\mathcal{D}(n)}
\def\rn{\mathcal{R}(n)}

\def\rnp{\mathcal{R}(n+1)}
\def\mo{\mathcal{O}}
\def\ma{\mathfrak{I}_0}
\def\mb{\mathfrak{I}_1}
\def\mc{\mathfrak{I}_2}

\begin{document}

\title[QUANTUM COHOMOLOGY AND TOTAL POSITIVITY]{QUANTUM COHOMOLOGY RINGS OF LAGRANGIAN AND ORTHOGONAL GRASSMANNIANS AND TOTAL POSITIVITY}
\author{Daewoong Cheong}
\address{Korea Institute for Advanced Study \\
207-43 Cheongryangri 2-dong\\
Seoul, 130-722, Korea} \email{daewoongc@kias.re.kr}
\date{October 26, 2006}

\begin{abstract}
We verify  in an elementary way a result of Peterson for the
maximal orthogonal and  Lagrangian Grassmannians, and then find
Vafa-Intriligator type formulas which compute their  $3$-point,
genus zero Gromov-Witten invariants. Finally we study the total
positivity of the related Peterson's varieties and show that
Rietsch's conjecture about the total positivity holds for these
cases.
\end{abstract}
\maketitle

\section{INTRODUCTION}

The (small) quantum cohomology ring $qH^*(X)$ of a projective
complex manifold $X$ is a certain deformation of the classical
cohomology ring of $X.$ The additive structure of the ring
$qH^*(X)$ is the same as that of the classical one, but the
multiplicative one is deformed from that of the classical
cohomology; the structure constants for the quantum multiplication
are the $3$-point Gromov-Witten invariants of genus $0,$ which
count equivalence classes of certain rational curves in $X.$

In 1997, D. Peterson announced in his MIT lectures that the
(complexified) quantum cohomology ring $qH_{\C}^*(G/P)$ of a
homogeneous manifold  $G/P$ is isomorphic to the coordinate ring
of a subvariety $\mathcal{Y}_{P}$ of  the so-called Peterson
variety $\mathcal{Y} \subset G^{\vee}/B^{\vee}$(see
\ref{peterson}). The proof of Peterson's general theorem remains
unpublished. In the general case it is not easy to see the quantum
cohomology rings concretely via those isomorphisms, because the
varieties $\mathcal{Y}_P$ are defined implicitly. However if $P$
is a minuscule parabolic subgroup, then we can replace the variety
$\mathcal{Y}_P$ by an isomorphic subvariety $\mathcal{V}_P $ of
$G^{\vee}$ whose coordinate ring is more explicit(see
\ref{peterson}).

In \cite{Riet1} Rietsch verified the result of Peterson for the
Grassmannian manifolds
 $GL_{n}(\C)/P$ in type $A$ by constructing in an elementary way an
isomorphism between the coordinate ring of $\mathcal{V}_{P}$ and
the quantum cohomology ring of $GL_n(\C)/P.$  In addition, using
this isomorphism, Rietsch derived for Grassmannian manifolds in
type $A$ the Vafa-Intriligator formula which gives the
Gromov-Witten invariants by evaluating certain symmetric functions
on some roots of unity. Then Rietsch identified the totally
nonnegative elements in the variety $\mathcal{V}_{P}$ and
characterized these elements via the nonnegativity of Schubert
basis elements on $\mathcal{V}_{P},$ which are determined by the
Schur polynomials.

\smallskip
 In this paper, we find analogues of Rietsch's results in
the other classical Lie types. Consider a complex vector space $V$
of dimension $N$ with a nondegenerate (symmetric or
skew-symmetric) bilinear form $Q.$ For $Q$ a symmetric bilinear
form and $N=2n+2$ , let $\oge$ be one component of the parameter
space of maximal isotropic subspaces of $V,$ and for $Q$ a
symmetric bilinear form and $N=2n+1,$ let $\ogo$ be  the parameter
space of maximal isotropic subspaces of $V.$ Recall that subspace
$W$ of $V$ is called isotropic if
 $Q(v,w)=0$ for all $v,w\in W.$ For  $N=2n$ and $Q$ a skew symmetric bilinear form, we denote by $\lg$
the parameter space of Lagrangian (i.e., maximal isotropic)
subspaces of $V.$ If we denote by $P_n$ all the parabolic
subgroups corresponding to `right end roots' in the Dynkin
diagrams for the classical algebraic groups with rank $n,$  then
$\ogo=SO_{2n+1}(\C)/P_n,$
 $\lg=Sp_{2n}(\C)/P_n$ and
 $\oge=SO_{2n+2}(\C)/P_{n+1}.$ It is well-known that  $\oge$
is isomorphic to $\ogo.$

\smallskip
 The cohomology ring $H^*(\lg)$ is generated by the Schubert classes $\sigma_\lambda,$ with $ \lambda$
 strict partitions $\lambda=(\lambda_1,...,\lambda_n)$(see \ref{subsec:sym poly}) such that
  $\lambda_1\leq n,$ and the quantum cohomology ring $\qhlg$ is
isomorphic, as $\Z[q]$-module, to $H^*(\lg)\otimes\Z[q].$
 The multiplication of $\qhlg$ is given by $$\sigma_\lambda \cdot \sigma_\mu=\sum_{\nu, d}
 <\sigma_\lambda,\sigma_\mu,\sigma_{\hat{\nu}}>_d \sigma_\nu q^d,$$
 where  the sum is taken over all strict partitions $\nu$ with $\nu_1\leq n$ and
 nonnegative integers $d,$ $\hat{\nu}$ denotes the partition whose parts complement the parts
 of $\nu$ in the  set $\{1,...,n\},$ so that the classes
 $\sigma_{\nu}$ and $\sigma_{\hat{\nu}}$ are Poincar\'{e} dual to
 each other, and $<\sigma_\lambda,\sigma_\mu,\sigma_{\hat{\nu}}>_d$ are Gromov-Witten
 invariants(See \ref{subsec:lg}).

\smallskip
 Similarly, the cohomology ring
$H^*(\oge)$ is generated by the Schubert classes $\tau_{\lambda}$
with $\lambda=(\lambda_1,...,\lambda_n)$ strict partitions such
that $\lambda_1\leq n$. The quantum cohomology ring $\qhoge$ is
isomorphic, as $\Z[q]$-module, to $H^*(\oge)\otimes\Z[q].$
 The multiplication of $\qhoge$ is given by $$\tau_\lambda \cdot \tau_\mu=\sum_{\nu, k}
 <\tau_\lambda,\tau_\mu,\tau_{\hat{\nu}}>_k\tau_\nu q^k,$$
 where the sum is taken over all strict partitions $\nu$ with $\nu_1 \leq n$ and nonnegative integers
 $k,$
and $<\tau_\lambda,\tau_\mu,\tau_{\hat{\nu}}>_k$ are
 Gromov-Witten invariants(See \ref{subsec:ogo}).

\smallskip
 To a triple $(G,G^{\vee},P)$ with $P$ a parabolic subgroup of $G,$
 we associate  a closed subvariety $\mathcal{V}_P$ of
 $G^{\vee}$(see \ref{peterson}).
  We write $\mathcal{V}_n^C,$ $\mathcal{V}_n^B$ and
$\mathcal{V}_{n+1}^D$ rather than $\mathcal{V}_P$  for the triples
$(\SOo, \Sp,P_n),$ $(\Sp, \SOo,P_n)$ and $(\SOe,\SOe,P_{n+1}),$
respectively. Let $\mo(\mathcal{V}_P)$ be the reduced coordinate
ring of $\mathcal{V}_P,$ and denote
 $qH^*_{\C}(X):=qH^*(X)\otimes \C.$
 The first main result of this paper is to  give the following isomorphisms explicitly:

\bn
\item$qH_{\C}^*(\ogo) \stackrel{\sim}{\rightarrow} \mathcal{O}[\mathcal{V}_{n}^C],$
\item $qH_{\C}^*(\lg) \stackrel{\sim}{\rightarrow} \mathcal{O}[\mathcal{V}_{n}^B],$
\item $qH_{\C}^*(\oge) \stackrel{\sim}{\rightarrow} \mo[\mathcal{V}_{n+1}^D].$
 \en

\smallskip

 As a second
result, we give Vafa-Intriligator type formulas to compute the
Gromov-Witten invariants for $\ogo$ $(\oge)$ and $\lg.$ \\
Let $\zeta=\zeta_n$ be the primitive $2n$-th root of unity, i.e.,
$\zeta_n=e^{\frac{\pi i}{n}},$ and $\mathcal{T}_n$ the set of all
$n$-tuples $J=(j_1,...,j_n),$ $-\frac{n-1}{2}\leq
j_1<\cdots<j_n\leq \frac{3n-1}{2}, $ such that
$\zeta^J:=(\zeta^{j_1},...,\zeta^{j_n})$ is an $n$-tuple of
distinct $2n$-th roots of $(-1)^{n+1}.$ Denote $\In:=\{ I\in
\mathcal{T}_n | \zeta^{i_k}\ne \zeta^{i_l}
\hspace{0.05in}\textrm{for} \hspace{0.05in}\textrm{all}
\hspace{0.05in} k,l=1,...,n\}$. For $J \in \mathcal{I}_n,$  define
$J^*$ to be a $n$-tuple $J^*=(j_1^*,...,j_n^*)\in \mathcal{I}_n$
such that the two sets $\{\zeta^{j_1},...,\zeta^{j_n}\}$ and
$\{\zeta^{n-j^*_1},...,\zeta^{n-j^*_n}\}$ enumerate all $2n$-th
roots of $(-1)^{n+1}.$  Let $\Qt_{\lambda}$ and $\Pt_\lambda$ be
the $\Qt$- and $\Pt$-polynomials of Pragacz-Ratajski(see
\ref{subsec:sym poly} or \cite{Pra and Rat 1}), and $P_\lambda$
the Hall-Littlewood symmetric polynomial.(\cite{Mac1}).

\bn
 \item
Vafa-Intriligator type formula for $\oge$($\cong \ogo$):\\
Given Schubert classes $ \tau_{\lambda},\tau_{\mu},$ $\tau_\nu$ of
$\oge,$ and a nonnegative integer $k,$  the Gromov-Witten
invariants $<\tau_{\lambda},\tau_{\mu},\tau_{\hat{\nu}}>_k$
 are given by
\begin{displaymath}
<\tau_{\lambda},\tau_{\mu},\tau_{\hat{\nu}}>_k=\frac{2^{l(\nu)+2k}}{(2n)^n}\sum_{m=0}^{a(\nu)}\sum_{J\in
\In}\widetilde{P}_{\lambda}(\zeta^{J})\widetilde{P}_{\mu}(\zeta^{J})P_{\nu(m)}(\zeta^{J^{*}})|
\textrm{Vand}(\zeta^J)|^2,
\end{displaymath}

whenever $|\lambda| + |\mu|=|\nu|+2nk$, and otherwise by
$<\tau_{\lambda},\tau_{\mu},\tau_{\hat{\nu}}>_k=0.$

Here  $a(\nu):=\lfloor \frac{n-l(\nu)} {2}\rfloor, $
$\nu(m):=((n)^{(2m)},\nu_1,...,\nu_l),$ and
$\textrm{Vand}(\zeta^J):=\prod_{k<l}(\zeta^{j_k}-\zeta^{j_l}).$

\bigskip
\item
Vafa-Intriligator type formula for $\lg:$ \\
Given Schubert classes $\sigma_{\lambda},\sigma_{\mu},$
$\sigma_{\nu}$ of $\lg,$ and a nonnegative integer $d,$ the
Gromov-Witten invariants
$<\sigma_{\lambda},\sigma_{\mu},\sigma_{\hat{\nu}}>_d$ are given
by
$$
<\sigma_{\lambda},\sigma_{\mu},\sigma_{\hat{\nu}}>_d=\frac{1}{2^d(2n+2)^{(n+1)}}\sum_{m=0}^{b(\nu)}
 \sum_{J\in
 \Im}\Qtl(\zeta^{J})\Qt_{\mu}(\zeta^{J})P_{\nu[m]}(\zeta^{J^{*}})|
 \textrm{Vand}(\zeta^J)|^2,
$$

whenever $|\lambda| + |\mu|= |\nu| +(n+1)d,$ otherwise
$<\sigma_{\lambda},\sigma_{\mu},\sigma_{\hat{\nu}}>_d=0.$

Here we denote $\zeta=\zeta_{n+1}.$ For the notation of $b(\nu)$
and $\nu[m],$ see the Proposition $\ref{fifth}.$\en

\smallskip
Lastly we explicitly describe the totally nonnegative parts of the
varieties $\mcvc$ and $\mcvb,$ and show that they are isomorphic
to $\R_{\geq0}.$ Furthermore for our cases we prove Rietsch's
conjecture, which states that it is only on the totally
nonnegative elements that the Schubert basis elements are
nonnegative.

 This paper is organized as follows. In Section $\ref{sec:2},$ we collect some basic facts about Lie
algebras in types $B,C,$ and $D,$ and give the definitions of
$\Qt$- and $\Pt$-polynomials of Pragacz and Ratajski. In Section
$\ref{sec:3},$ we explain the quantum cohomology rings of
Lagrangian and orthogonal Grassmannians, and describe the result
of Peterson. In Section $\ref{sec:4},$ we give
 the Bruhat decomposition of the stabilizer $(U^{+})^e$ of $e$ in
$U^{+},$ determine by which regular functions the closure of a
cell is defined, and we state the main theorems on Peterson's
result, theorems \ref{thm;ogo},  \ref{2nd} and \ref{thm:oge}, and
partially prove them. In Section $\ref{sec:5},$ we
 characterize elements of $\mcvb,$ $\mcvc$ and
$\mcvd,$ and describe how subvarieties of these varieties are
positioned.
 In Section $\ref{sec:6},$ we prove some
orthogonality formulas for $\Qt$- and $\Pt$-polynomials at roots
of unity, use these formulas to complete the theorems
\ref{thm;ogo}, \ref{2nd} and $\ref{thm:oge}.$ In Section
$\ref{sec:7},$  we derive the Vafa-Intriligator type formulas for
Lagrangian and orthogonal Grassmannians, and also we give an
analogue of Poincar\'{e} duality pairing on $qH^*(\oge)$ and
$qH^*(\lg).$ In Section $\ref{sec:8},$ we give a quick review of
total positivity of $U^+$ and describe the totally positive parts
of $\mcvb$ and $\mcvc.$

\section{PRELIMINARIES}\label{sec:2}
\subsection{Notations}
Let $G$ be a complex semisimple algebraic group with rank $n$. Let
$ B=B^{-}$ and $B^{+}$ be opposite Borel subgroups and $U^{-}$ and
$U^{+}$ be the unipotent radicals of $B^{-}$ and $B^{+}$
respectively. Let $T=B^{-}\cap B^{+}$ be a maximal torus of $G.$
Let $\mfg,$   $\mfbm,$ $\mfb,$ $\mfum$, $\mfu$ and $\mft$ be the
Lie algebras of $G,$ $B^{-},$ $B^{+},$ $U^{-},$ $U^{+},$ and $T$
respectively, so that the Lie algebra $\mfg$ has the Cartan
decomposition $\mfg=\mfum \oplus \mft \oplus \mfu.$  Let $e_i,$
$h_i$ and $f_i,$ $i=1,...,n,$ be the standard generators of
$\mfg,$ and let $A=(a_{ij})$ be the Cartan matrix of $\mfg.$

We denote by  $\triangle$  the root system for $\mfg$ and by
$\triangle^+$  the set of positive roots in $\triangle$ and  by
$\Pi:=\{\alpha_1,...,\alpha_n\} \subset \triangle$ the set of
simple roots , which are defined by $\alpha_j({h}_i)=a_{ij}.$ Let
$G^{\vee}$ be the Langlands dual of $G,$ and denote the
counterparts in $G^{\vee}$ of  the above subgroups and their Lie
algebras by the same alphabets as above with $^{\vee}.$ The $
Weyl$ $group$ $W$ of $G$ is defined as $W=\textrm{Norm}_G(T)/T$.
The action of $W$ on T induces the action of $W$ on $\mft$ and the
dual action on $\mft^{*}$. If $W$ is generated by simple
reflections $s_i$, $i\in I$, the action of $W$ on $\mft^{*}$ is
given by $s_i(\gamma)=\gamma-\gamma({h}_i)\alpha_i$, where
$\gamma\in \mft^{*}$, ${h}_i \in \mft$, and $\alpha_{i}\in \Pi.$
Note that the Weyl groups of  $G$ and $G^{\vee}$ are identical to
each other.

Let $w_0$ be the longest element of $W,$ and for a parabolic
subgroup $P,$ let $W_P$ be the Weyl group of $P$ and $w_P$ the
longest element of $W_P.$ Then $w^P:=w_0 w_P$ is the  minimal
length representative in the coset $w_0W_P\in W/W_P.$ If a
parabolic subgroup $P$ corresponds to a subset $J\subset
I:=\{1,...,n\},$ we write $P=P_J,$ and we use alternatively $J$
and $P$ in subscript or superscript of the above notations, e.g.,
$ w^J=w^P.$ If $P$ is a maximal parabolic subgroup corresponding
to a fundamental weight $\kappa_r$ for some $r\in I,$ i.e.,
$P=P_J$ for $J=I\setminus \{r\},$ then we write $P_r$ for $P_J.$

 A fundamental weight $\kappa_r$ is $minuscule$ if  each weight of the fundamental representation
 $V_r$ corresponding to $\kappa_r$ is extremal, that is, if it  is of the form $w\cdotp \kappa_r$
 for some $w\in W.$
 A maximal parabolic subgroup $P_r$ is called $minuscule$ if the fundamental weight
  $\kappa_r$ is minuscule.
 We have the following list of the minuscule weights in classical
 Lie groups.\\
 \indent \indent \indent \indent Type $A_n:$ Every fundamental weight is minuscule.\\
 \indent \indent \indent \indent Type $B_n:$ $\kappa_n$ is minuscule.\\
 \indent \indent \indent \indent Type $C_n:$ $\kappa_1$ is minuscule.\\
\indent  \indent \indent \indent Type $D_n:$
$\kappa_1,\kappa_{n-1},\kappa_{n}$ are minuscule.\\
 From the above list, we note that the parabolic subgroups $P_n \subset \SOo$ and
 $P_{n+1} \subset \SOe$ are minuscule, but the parabolic subgroup $P_n \subset
 Sp_{2n}(\C)$ is not minuscule. We will mainly work with these
 three parabolic subgroups.

\smallskip
 In the following three subsections, we summarize some basic facts about
 Lie algebras in type $B,C$ and $D,$ which will be used in later
 sections. We will take other ordered basis for $V$ than usual Lie theory
 literatures do, in order to embed the Lie groups in $B,C$ and $D$ into
 the Lie group in type $A$ more $naturally$.

\subsection{ $Sp_{2n}(\C)$ and $\mathfrak{sp}_{2n}(\C)$}~\label{basic:sp}
Let $V$ be a $2n$-dimensional complex vector space equipped with a
nondegenerate skew-symmetric bilinear form Q. Then the symplectic
Lie group $Sp_{2n}(\C)$ is defined to be the group of
automorphisms $A$ of $V$ preserving $Q$, i.e., $Q(Av,Aw)=Q(v,w)$
for all $v$, $w\in V$. The Lie algebra $\mathfrak{sp}_{2n}(\C)$ of
$Sp_{2n}(\C)$ is the vector space of endomorphisms $X:
V\rightarrow V$ such that $Q(Xv,w)+Q(v,Xw)=0$ for all $v$, $w\in
V$. Choose  an ordered  basis $\xi_1,...,\xi_{2n}$ for $V$ such
that

$$ Q(\xi_{i},\xi_{i^\ast})=-Q(\xi_{i^\ast},\xi_i)=(-1)^{i+1},$$
 and
$$Q(\xi_i,\xi_j)=0 \hspace{0.07in}\textrm{if}\hspace{0.07in} i+j\neq 2n+1,$$
where $i^\ast=2n+1-i.$ Let $J$ be the $(2n\times 2n)$-matrix
defined by \be \label{sym:J} J_{i,j}:=
(-1)^{i+1}\delta_{i,j^\ast}, \hspace{0.05in} i=1,...,2n
 \ee
Then $Sp_{2n}(\C)$ is the group of $2n\times 2n$ matrices  $A$
with $J=A^t\cdotp J\cdotp A$, and $\mathfrak{sp}_{2n}(\C)$ is the
space of  $2n\times 2n$ matrices $X$ with $X^t\cdotp J+J\cdotp
X=0$.

For positive integers $i,j,$ let $E_{i,j}$ be the $2n\times
2n$-matrices such that $(i,j)$-th entry is $1$, and $0$ elsewhere,
and let $\bar{h}_i:=E_{i,i}-E_{i^*,i^*}.$ As Cartan subalgebra
$\mathfrak{t}$ for $\mathfrak{sp}_{2n}(\C),$ we take the
subalgebra generated by $\bar{h}_i,$ $i=1,...,n.$ Let $l_1,...,l_n
$ be a basis of $ \mathfrak{t}^{*}$ dual to
$\bar{h}_1,...,\bar{h}_n.$ The set $\triangle=\{\pm l_{i} \pm
l_{j} \mid 1\leq i, j\leq n \}$ forms a root system for
$\mathfrak{sp}_{2n}(\C).$ The set $\Pi$ of simple roots consists
of $\alpha_i:=l_i-l_{i+1}$ for $i=1,...,n-1$, and
$\alpha_n:=2l_n,$ and the simple root vectors (\cite{BZ1}) are
$$e_i=E_{i,i+1}+E_{i^*-1,i^*} \hspace{0.07in}\textrm{for}\hspace{0.07in}i=1,..n-1,$$
$$e_n=E_{n,n+1}.$$  The fundamental dominant
weights corresponding to $\alpha_i$ are $\kappa_i=l_1+\cdots+l_i$,
$i=1,...,n,$ and the fundamental representation $V_i$
corresponding to $\kappa_i$ is the subspace of $\wedge^i V$
generated by a highest weight vector $\xi_1\wedge \cdots \wedge
\xi_i.$

The Weyl group $W_n$ for the root system $C_n$ is an extension of
the symmetric group\\
$S_n=<s_1,...,s_{n-1}>$ by $\dot{s}_n,$ which acts on the right,
in the notation of bar permutations,
$$(u_1,...,u_n)\dot{s}_n=(u_1,...,u_{n-1},\bar{u}_n).$$
The generators $s_1,...,s_{n-1},$ and $\dot{s}_n$  act naturally
on $\mft^{\ast}$; for $i\leq n-1$, $s_i$ interchanges $l_i$ and
$l_{i+1}$, and $\dot{s}_n$ interchanges $\pm l_n$ and $\mp l_n$.
Note that the maximal length element $w_0$ in $W_n$ is
$(\bar{1},\bar{2},...,\bar{n}).$

\subsection{$SO_{2n+1}(\C)$ and
$\mathfrak{so}_{2n+1}(\C)$}\label{odd orthogonal group}
\label{subsec 2n+1} Let $V$ be a $(2n+1)$-dimensional complex
vector space equipped with a nondegenerate symmetric bilinear form
$Q: V \times V\rightarrow \C$. The orthogonal group
$SO_{2n+1}(\C)$ is defined to be the group of automorphisms $A$ of
$V$ of determinant $1$ preserving $Q$, that is, $Q(Av,Aw)=Q(v,w)$
for all $v,w \in V$, and its Lie algebra
$\mathfrak{so}_{2n+1}(\C)$ is a vector space of endomorphisms
$X:V\rightarrow V$ such that $Q(Xv,w)+Q(v,Xw)=0$ for all $v,w \in
V$.  Let $\xi_1,...,\xi_{2n+1}$ be a basis for $V$ such that
$$Q(\xi_i,\xi_{i^\ast})=Q(\xi_{i^\ast},\xi_i)=(-1)^{i+1},$$ and
$$Q(\xi_i,\xi_j)=0 \hspace{0.07in}\textrm{if}\hspace{0.07in} i+j\neq 2n+2,$$
where $i^\ast=2n+2-i.$ Let $J$ be the $(2n+1)\times (2n+1)$-matrix
defined by
$$J_{i,j}:=(-1)^{i+1}\delta_{i,j^\ast},\hspace{0.06in}i,j=1,...,2n+1.$$
Then $SO_{2n+1}(\C)$ is the group of matrices $A$ satisfying the
relations $J=A^t\cdotp J\cdotp A,$ and the Lie algebra
$\mathfrak{so}_{2n+1}(\C)$ is the space of matrices X satisfying
the relation $X^t\cdotp J+J\cdotp X=0 $. Define $E_{i,j},$ $l_i$
and $\bar{h}_i$ as in \ref{basic:sp}. As Cartan subalgebra we take
the subalgebra generated by $\bar{h}_i,$ $i=1,...,n.$  The set
$\triangle=\{\pm l_i \pm l_j \mid i,j=1,...,n, \hs \textrm{and}
\hspace{0.08in}i\neq j \}\cup \{\pm l_i \mid i=1,...,n\}$ forms a
root system for $\mathfrak{so}_{2n+1}(\C).$ The set of simple
roots consists of $\alpha_i:=l_i-l_{i+1}$ for $i=1,...,n-1,$ and
$\alpha_n:=l_n.$ The simple root vectors (\cite{BZ1}) are
$$e_i=E_{i,i+1}+E_{i^*-1,i^*} \hspace{0.07in}
\textrm{for}\hspace{0.07in} i=1,...n-1,$$
$$e_n=\sqrt{2}(E_{n,n+1}+E_{n+1,n+2}).$$  The
fundamental dominant weights are $\kappa_i=l_1+ \cdots +l_i$ for
$i=1,...,n-1$, and $\kappa_n=\frac{1}{2}(l_1+ \cdots + l_n)$. To
each $\kappa_i,$ $i=1,...,n-1,$ there corresponds fundamental
representation $V_i$ of $\mathfrak{so}_{2n+1}(\C)$ which is  the
subspace of $\wedge^i V$ generated by a highest weight vector
$\xi_1\wedge\cdots\wedge \xi_i$, and to $\kappa_n$  the spin
representation $V_n^s$(\cite{FH}). The group $SO_{2n+1}(\C)$ has
$V_1,...,V_{n-1}$ as fundamental representations, but not $V_n^s.$
But there is a representation $V_n$ which corresponds to the
weight $l_1+\cdots + l_n.$  This is the subspace of $\wedge^n V$
generated by a highest weight vector $\xi_1\wedge\cdots \wedge
\xi_n.$ We consider $V_n$ as fundamental representation as well.

As a matter of convenience for later use, we give basic facts on
$SO_{2n+2}(\C)$ and $\mathfrak{so}_{2n+2}(\C)$ rather than
$SO_{2n}(\C)$ and $\mathfrak{so}_{2n}(\C).$
\subsection{ $SO_{2n+2}(\C)$ and $\mathfrak{so}_{2n+2}(\C)$}\label{sec:basic D}

Let  $V$ be a $(2n+2)$-dimensional vector space with a
nondegenerate bilinear form $Q.$ The definitions of
$SO_{2n+2}(\C)$ and $\mathfrak{so}_{2n+2}(\C)$ are the same as
those of $SO_{2n+1}(\C)$ and $\mathfrak{so}_{2n+1}(\C).$ We take
an ordered basis $\xi_1,...,\xi_{2n+2}$ for $V$ such that
$$Q(\xi_i, \xi_{i^*})=Q(\xi_{i^\ast},\xi_i)=1,$$ and
$$Q(\xi_i,\xi_j)=0\hspace{0.06in}\textrm{if}\hspace{0.06in}i+j\ne 2n+3,$$
where $i^\ast=2n+3-i.$ Let $E_{i,j},$ $\bar{h}_i$ and $l_i$ be
defined as in \ref{basic:sp}. As Cartan subalgebra we take  the
subalgebra generated by $\bar{h}_i,$ $i=1,...,n+1.$ The set of
roots of $\mathfrak{so}_{2n+2}(\C)$ is
$$\triangle=\{\pm l_i\pm l_{j}\mid i,j=1,...,n+1,i\neq j\} $$  and the set $\Pi$ of simple roots
consists of $\alpha_i=l_i-l_{i+1},$ $i=1,...,n,$ and
$\alpha_{n+1}=l_{n}+l_{n+1}.$  The corresponding simple root
vectors are
$$e_i=E_{i,i+1}-E_{i^{*}-1,i^{*}}\hspace{0.06in}\textrm{for}\hspace{0.07in} i=1,...,n,$$
$$e_{n+1}=E_{n,n+2}-E_{n+1,n+3}.$$
Denote $\omega_i:=l_1+\cdots+l_i$ for $i=1,...,n,$ and
$\omega_{n+1}^{\pm}:=l_1+\cdots+l_{n}\pm l_{n+1}.$ Then the
fundamental dominant weights $\kappa_i$  corresponding to
$\alpha_i$ are as follows: $\kappa_i=\omega_i$ for $i=1,...,n-1,$
$\kappa_{n}=\frac{1}{2}\omega_{n+1}^{-},$ and
$\kappa_{n+1}=\frac{1}{2}\omega_{n+1}^{+}.$ To a fundamental
weight $\kappa_i,$ $i=1,...,n-1,$ there corresponds the
fundamental representation $V_i$ of $\mathfrak{so}_{2n+2}(\C)$
which is the subspace of $\wedge^iV$ generated by a highest weight
vector $\xi_1\wedge\cdots\wedge \xi_i,$  and to the fundamental
weights $\kappa_{n}$ and $\kappa_{n+1}$ there correspond the spin
representations $V_{n}^s$ and $V_{n+1}^{s},$
respectively(\cite{FH}). As fundamental representations, the group
$SO_{2n+2}(\C)$ has $V_i,$ $i=1,...,n-1,$ but not the spin
representations $V_{n}^s$ and $V_{n+1}^s.$ But there are
representations $V_{n}$ and $V_{n+1}^{\pm}$ of $SO_{2n+2}(\C),$
which correspond to the weights $\omega_{n}$ and
$\omega_{n+1}^{\pm},$ respectively. The representation $V_{n}$ is
the subspace of $\wedge^{n}V$ generated by a highest weight vector
$\xi_1\wedge \cdots \wedge \xi_{n},$ and $V_{n+1}^{\pm}$ is the
subspaces of $\wedge^{n+1}V$ generated by highest weight vectors
$\xi_1 \wedge \cdots \wedge \xi_{n+1} $ for $(+)$ sign   and
$\xi_1\wedge \cdots \wedge \xi_{n}\wedge \xi_{n+2}$ for $(-)$
sign. Note that these representations generate all the
representations of $SO_{2n+2}(\C), $ with the following relation (
p.379 in \cite{FH}). \be \label{b.1} (V_{n+1} ^{+} \oplus V_{n-1}
\oplus V_{n-3} \oplus \cdots ) \otimes (V_{n+1}^{-} \oplus V_{n-1}
\oplus V_{n-3}\oplus \cdots )=(V_{n}\oplus V_{n-2}\oplus
\cdots)^2.\ee We will refer to these $(n+2)$ representations as
fundamental ones of $SO_{2n+2}(\C).$

The Wely group $\widetilde{W}_{n+1}$ for type $D_{n+1}$ is an
extension of the symmetric group $S_{n+1}=<s_1,...,s_{n}>$ by an
element $\tilde{s}_{n+1}$ which acts on the right by
$$(u_1,...,u_{n+1})\tilde{s}_{n+1}=(u_1,...,u_{n-1},\bar{u}_{n+1},\bar{u}_{n}).$$
 Note that  the maximal
length element $\widetilde{w}_0$ of $\widetilde{W}_{n+1}$ is given
by \be\widetilde{w}_0= ~\label{w0}
 \left
\{\begin{array}{cc}
 (\bar{1},...,\overline{n+1})& \hspace{0.08in}\textrm{if n is odd},\\
 (\bar{1},...,\bar{n},n+1) & \hspace{0.08in}\textrm{if n is
 even}.
 \end{array}\right. \ee
 The action of $\widetilde{W}_{n+1}$ on
$\mathfrak{t}^*$ is given as follows; for $i=1,...,n,$ $s_i$
interchanges $l_i$ and $l_{i+1},$ and $\tilde{s}_{n+1}$
interchanges $ \pm l_{n}$ and $\mp l_{n+1}.$

\subsection{Symmetric Polynomials}\label{subsec:sym poly}
A $partition$ $\lambda$ is a weakly decreasing sequence
$\lambda=(\lambda_1,\lambda_2,...,\lambda_m)$ of nonnegative
integers. A $Young \hs diagram $ is a collection of boxes,
arranged in left-justified rows, with a weakly decreasing number
of boxes in each row. To a partition
$\lambda=(\lambda_1,...,\lambda_m),$ we associate a Young diagram
whose $i$-th row has $\lambda_i$ boxes. The nonzero $\lambda_i$ in
$\lambda=(\lambda_1,...,\lambda_m)$ are called the $parts$ of
$\lambda.$ The number of the parts of $\lambda$ is called the
$length$ of $\lambda,$ denoted by $l(\lambda)$; the sum of the
parts of $\lambda$ is called the $weight$ of $\lambda$, denoted by
$|\lambda|.$  For positive integers $m$ and  $n,$ denote by
$\mathcal{R}(m,n)$ the set of all partitions whose Young diagram
fits inside an $m\times n$ diagram, which is the Young diagram of
the partition $(n^m).$ A partition
$\lambda=(\lambda_1,...,\lambda_m)\in \mathcal{R}(m,n)$ is called
strict if $\lambda_1>\cdots>\lambda_m.$ Denote by
$\mathcal{D}(m,n)$ the set of all strict partitions in
$\mathcal{R}(m,n).$ If $m=n$, then we write $\mathcal{R}(n)$ and
$\mathcal{D}(n)$ for $\mathcal{R}(n,n)$ and $\mathcal{D}(n,n)$,
respectively. If $\lambda \in \mathcal {D}(n)$, denote by
$\widehat{\lambda}$ the partition whose parts complements the
parts of $\lambda$ in the set $\{1,...,n\}$.  We define $\Qt$- and
$\Pt$-polynomials of Pragacz and Ratajski  as follows. For
$X_n:=(x_1,...,x_n),$  set $\Qt_i(X_n):=E_i(X_n)$, the $i$-th
elementary symmetric function. Given two integers $i$ and $ j$
with $i\geq j$, define
\begin{displaymath}
\Qt_{i,j}(X_n)=\Qt_i(X_n)\Qt_j(X_n)+2\sum_{k=1}^{j}(-1)^k
\Qt_{i+k}(X_n)\Qt_{j-k}(X_n).
\end{displaymath}

 Finally, for any partition
$\lambda$ of length $l=l(\lambda)$, not necessarily strict, define
\begin{displaymath}
\Qt_{\lambda}(X_n)=\textrm{Pfaffian}[\Qt_{\lambda_i,\lambda_j}(X_n)]_{1\leq
i,j \leq n},
\end{displaymath}
 where $n=2\lfloor(l+1)/2\rfloor.$

The $\Qt$-polynomials satisfy the following properties which will
be used in later sections.
\begin{enumerate}
\item $\Qt_{i,i}(X_n)=E_i(x_1^2,...,x_n^2).$ \item For partitions
$\lambda=(\lambda_1,...,\lambda_l)$ and
$\lambda^{\prime}=(\lambda_1,\lambda_2,...,j,j,...,\lambda_l)=\lambda
\cup (j,j),$
$$\Qt_{\lambda^{\prime}}(X_n)=\Qt_{j,j}(X_n)\Qt_{\lambda}(X_n). $$\label{fac:Q}
\item
For any $\lambda \in \mathcal{D}(n),$
$$\Qtl(X_n)\Qt_n(X_n)=\Qt_{(n,\lambda_1,...,\lambda_l)}(X_n).$$
\label{fac:Qspe}
\end{enumerate}

Given $\lambda$, not necessarily strict,  $\Ptl$ is defined by
$$\Ptl(X_n):=2^{-l(\lambda)}\Qtl(X_n).$$
Note that $\Pt$-polynomials enjoy the factorization properties
$(2)$ and $(3)$ of $\Qt.$ For $i=1,...,n,$ let $H_i(X_n)$ be the
$i$-th complete symmetric function. Then for any partition
$\lambda,$ the Schur polynomial $S_{\lambda}(X_n)$ is defined by
$$S_{\lambda}(X_n):=\textrm{Det}[H_{\lambda_i-j+i}(X_n)]_{1\leq i,j
\leq n },$$ where $H_0(X_n)=1$ and $H_{k}(X_n)=0$ for $k<0.$ See
\cite{Pra and Rat 1} for a further reference on $\Qt$- and
$\Pt$-polynomials, and \cite{Mac1} for the Schur polynomials.

\section{QUANTUM COHOMOLOGY RINGS AND PETERSON'S RESULT}\label{sec:3}
 In this section, we define the quantum cohomology rings of $\lg$
and $\oge,$  give presentations of these rings and describe
Peterson's result on the quantum cohomology ring of $G/P.$
\subsection{ Quantum Cohomology ring of Lagrangian Grassmannian} \label{subsec:lg}
 To describe the quantum cohomology of $\lg$, we define the
Schubert varieties indexed by the strict partitions in $\dn.$ Let
$V$ be a complex vector space of dimension $2n$ with a
nondegenerate skew-symmetric form. Let $ F_{\mbox{\boldmath{.}}}$
be a fixed complete isotropic flag of subspaces $F_i$ of $V$:
$$\Fdot : 0=F_0\subset F_1\subset \cdots \subset F_n\subset V,$$
where dim$(F_i)=i$  for each $i,$ and $F_n$ is Lagrangian. Define
the Schubert variety $X_{\lambda}(\Fdot)$ as the locus of $ \Sigma
\in \lg$ such that \be \label{schubert} \textrm{dim}(\Sigma \cap
F_{n+1-\lambda_i})\geq i \trm i=1,...,l(\lambda).\ee Then
$X_{\lambda}(\Fdot)$ is a subvariety of $\lg$ of complex
codimension $| \lambda |.$ Let $\sigma_{\lambda}$ be the class of
$X_{\lambda}(\Fdot)$ in the cohomology group
$H^{2|\lambda|}(\lg).$ Then $\sigma_{\lambda}$ with $\lambda\in
\dn$ form an additive basis for $H^{*}(\lg).$

A rational map of degree $d$ to $\lg$ is a morphism
$f:\mathbb{P}^1\rightarrow \lg$ such that
$$\int_{\lg}f_*[\mathbb{P}^1]\cdotp \sigma_1=d.$$ Given an integer $d\geq0$ and partitions
$\lambda, \mu,$  $\nu\in \mathcal{D}(n)$ with $ |\lambda| + |\mu|
+ |\nu| = \textrm{dim} (\lg)+d(n+1),$ the Gromov-Witten invariants
$<\sigma_{\lambda},\sigma_{\mu},\sigma_{\nu}>_d$  are defined as
the number of rational maps $f:\mathbb{P}\rightarrow \lg$ such
that $f(0)\in X_{\lambda}(\Fdot),$ $f(1)\in X_{\mu}(
G_{\mbox{\boldmath{.}}}),$ and $f(\infty)\in X_{\nu}(
H_{\mbox{\boldmath{.}}})$, for given isotropic flags $\Fdot$,
$G_{\mbox{\boldmath{.}}},$ and $H_{\mbox{\boldmath{.}}}$ in
general position. The quantum cohomology ring $qH^{*}(\lg)$ is
isomorphic
 to $H^{*}(\lg)\otimes_{\Z} \Z[q]$
as a $\Z[q]$-module, where $q$ is a formal variable of degree
$(n+1).$ The multiplication in $qH^{*}(\lg)$ is given by the
relation \be~\label{multi} \sigma_\lambda \cdotp
\sigma_\mu=\sum<\sigma_\lambda,\sigma_\mu,\sigma_{\hat{\nu}}>_d
\sigma_\nu q^d,\ee  where the sum is taken over $d\geq0$ and
partitions $\nu$ with $|\nu|=|\lambda|+|\mu|-(n+1)d.$ We  have the
following presentation of the quantum cohomology ring of $\lg,$
due to Kresch and Tamvakis.
\begin {theorem}[\cite{KT2}] \label{quantum coho:lg}
The ring $qH^{*}(\lg)$ is presented as a quotient of the
polynomial ring $\Z[\sigma_{1},...,\sigma_n, q]$ by the relations
$$\sigma_r^2+2\sum_{i=1}^{n-r}(-1)^i\sigma_{r+i}\sigma_{r-i}=(-1)^{n-r}\sigma_{2r-n-1}q$$
for $1\leq r \leq n.$
\end{theorem}

See \cite {KT2} for more details on the quantum cohomology ring of
$\lg.$

\subsection{Quantum cohomology ring of Orthogonal  Grassmannian}\label{subsec:ogo}  Let $V$ be a complex vector space of
dimension $2n+2$ with a nondegenerate symmetric form. The Schubert
varieties $X_{\lambda}(\Fdot)$ are parametrized by partitions
$\lambda \in \mathcal{D}(n)$, and are defined by the same equation
$(\ref{schubert})$ as before, relative to an isotropic flag $
F_{\mbox{\boldmath{.}}}$ in $V.$  Let $\tau_\lambda$ be the
cohomology class of $X_{\lambda}(\Fdot)$. Then $\tau_\lambda$ for
$\lambda \in \mathcal{D}(n)$ form a $\Z$-basis for $H^*(\oge).$
The cohomology ring $H^*(\oge)$ can be presented in terms of
$\Pt$- polynomials. More precisely, let $\Lambda_n$ denote
$\Z$-algebra generated by the polynomials $\Pt_{i}(X_n)$ for all
$i=1,...,n$. Then the map from $\Lambda_n$ to $H^{*}(\oge)$
sending $\Ptl(X_n)$ to $\tau_{\lambda}$ is a surjective ring
homomorphism with the kernel generated by the polynomials
$\Pt_{i,i}(X_n)$ for all $i=1,...,n$ ([17, Sect.6] and \cite{Pra
and Rat 1}).

\smallskip
 In this case, the Gromov-Witten invariants
$<\tau_\lambda, \tau_\mu, \tau_\nu>_k$ are defined for
$|\lambda|+|\mu|+|\nu|=\textrm{deg}(\oge)+2nk$ and count the
number of rational maps $f:\mathbb{P}^1\rightarrow \oge $ of
degree $k$ such that $f(0)\in X_{\lambda}(\Fdot),$ $f(1)\in
X_{\mu}( G_{\mbox{\boldmath{.}}}),$ and $f(\infty)\in X_{\nu}(
H_{\mbox{\boldmath{.}}})$, for given isotropic flags $\Fdot$,
$G_{\mbox{\boldmath{.}}},$ and $H_{\mbox{\boldmath{.}}}$ in
general position. The quantum cohomology ring of $\oge$ is
isomorphic to $H^*(\oge)\otimes \Z[q]$ as a $\Z[q]$-module. The
multiplication in $qH^*(\oge)$ is given by the relation
\be~\label{multi:2} \tau_\lambda \cdotp
\tau_\mu=\sum<\tau_\lambda,\tau_\mu,\tau_{\hat{\nu}}>_k\tau_\nu
q^k,\ee where the sum is taken over $k \geq0$ and partitions $\nu$
with $|\nu|=|\lambda|+|\mu|-2nk.$

\begin{theorem}[\cite{KT1}]\label{quantum coho:oge}
The quantum cohomology ring $qH^{*}(\oge)$ is presented as a
quotient of the polynomial ring $\Z[\tau_1,...,\tau_n,q]$ modulo
the relations $\tau_{r,r}=0$ for $r=1,...,n-1$ together with the
quantum relation $\tau_n^2=q$, where \be \label{tau_ii}
\tau_{r,r}:=\tau_r^2+2\sum_{i=1}^{r-1} (-1)^i \tau_{r+i}
\tau_{r-i} + (-1)^r \tau_{2r}.\ee
\end{theorem}

See \cite{KT1} for more details on the quantum cohomology of
$\oge$ and $\ogo.$

\subsection{Peterson's results of quantum
cohomology ring of $G/P$.}\label{peterson}

For $i\in I,$ let $e_i^\vee $ be weight vectors of simple roots
for $\mathfrak{g}^\vee,$  ${(e_i^{\vee}})^{*} $ the linear
functional which is one on $e_i^{\vee}$ and zero on the weight
vectors of all the other roots, and set $$e^\vee:=\sum_{i\in
I}e_i^\vee \in \mathfrak{g}^\vee, \hspace{0.07in}\textrm{and}
\hspace{0.07in}(e^\vee)^*:=\sum_{i\in I}(e_i^{\vee})^* \in
({\mfg^{\vee}})^{*}.
$$
For $g\in G^\vee,$  let $Ad^{*}_g$ denote the coadjoint action of
$g$ on  $ (\mfg^{\vee})^{*},$ and let $((U^\vee)^+)^{e^\vee}$ be
the stabilizer of $e^\vee$ in $(U^\vee)^+.$ Then the Peterson
variety is defined as
$$\mathcal{Y}:=\{g(B^{-})^{\vee} \in G^{\vee}/(B^{\vee})^{-} \mid
Ad_g^{*}(e^\vee)^* \hspace{0.07in} \textrm{vanishes
on}\hspace{0.07in}
[(\mathfrak{u}^{\vee})^{-},(\mathfrak{u}^{\vee})^{-}]\}.$$ For a
parabolic subgroup $P\subseteq G,$  define $\mathcal{Y}_P$ as
$$\mathcal{Y}_P:=\mathcal{Y}\times_{G^{\vee}/B^{\vee}}
[(B^{\vee})^{+}w_{P}(B^\vee)^-/(B^\vee)^-] \hspace
{0.08in}(\textrm{scheme-theoretic intersection}).$$
 Then the subvarieties
$\mathcal{Y}_P$ form strata of the variety $\mathcal{Y}$, i.e,
$$\mathcal{Y}=\bigsqcup_{P}\mathcal{Y}_P,$$
where $P$ runs over the set of all parabolic subgroups $P$ of $G$
containing $B,$ and it is known that the variety $\mathcal{Y}_P$
need not be reduced, but is a local complete intersection. One of
the key results of Peterson's may be stated as follows:

\smallskip

The quantum cohomology ring $qH_\C^{*}(G/P)$ of a homogeneous
space $G/P$ is isomorphic to the  coordinate ring  of
$\mathcal{Y}_P.$

\smallskip If $P$ is a minuscule parabolic subgroup of $G,$
Peterson's result goes further. For a parabolic subgroup $P$ of
$G,$ not necessarily minuscule, define
$$\mathcal{V}_P:=((U^\vee)^+)^{e^\vee}\cap \overline{(B^\vee)^-w^P
(B^\vee)^-}.$$

\smallskip
If $P$ is minuscule, the quantum cohomology ring $qH_\C^{*}(G/P)$
of a homogeneous space $G/P$ is isomorphic to the coordinate ring
of $\mathcal{V}_P.$

\smallskip
This follows from an isomorphism of two varieties $\mathcal{V}_P
\stackrel{\sim}{\rightarrow} \mathcal{Y}_P,$ $u \mapsto
uw_{P}(B^\vee)^{-}$ (\cite{Pete2}).

\smallskip This result was verified in an elementary way by Rietsch
for minuscule parabolic subgroups  and  general parabolic
subgroups of $SL_N(\C)$ in \cite{Riet1} and \cite{Riet2},
respectively. See\cite{Pete1} and \cite{Riet3} for more details on
Peterson's results.

\section{COMPARING TWO PRESENTATIONS OF QUANTUM COHOMOLOGY RINGS}\label{sec:4}
In this section we will show  that two presentations for $\ogo,
\oge$ and $\lg$ by A. Kresch and H. Tamvakis and by D. Peterson
are equivalent to each other. This provides an independent proof
of Peterson's results for Lagrangian and orthogonal Grassmannians.

\subsection{$(U^{+})^e$ and its intersection with Bruhat
cells.}\label{Bruhat cells 1}

Throughout  this paper, we fix the principal nilpotent element
$e:=\sum_{i\in I}e_i\in \mfg$. Let $(U^+)^e:=\{u\in (U^+)^e\mid
ueu^{-1}=e$\}, the stabilizer of $e$ in $U^+$. This is an abelian
subgroup of $U^+$ of dimension equal to the rank of $G$
(\cite{Ko}).

\smallskip
\begin{lemma}[\cite{Riet1}]~\label{lem:1}
The elements $w_J$ can be characterized by  $\{w_J \in W|\hs
J\subseteq I \}=\{w\in W | \hs w\cdotp \Pi \subseteq(-\Pi)\cup
\triangle^{+}\}$.
\end{lemma}
\bp For the proof, we refer to Lemma $2.2$ in \cite{Riet1} .\ep
\smallskip
\begin{lemma}[\cite{Riet1}]~\label{lem:2}
Bruhat decomposition induces
 $(U^{+})^e$=$\bigsqcup_{J\subseteq I}(U^{+})^e \cap
B^{-}w^{J}B^{-}.$
\end{lemma}

\bp If $u\in (U^+)^e,$ then $u\in B^{-}w_0 wB^{-}$ for some $w\in
W.$ Write $u=b_1w_0wb_2$ for some $b_1\in B^{-}$, and $b_2\in
U^{-}.$ Since $ueu^{-1}=e,$ i.e., $u \cdotp e=e,$  we have
$$wb_2\cdotp e =w_0b_1^{-1}\cdotp e.$$ But note that $b_2\cdotp e
=e+x$ for some $x\in \mfbm,$  and $w_0b_1^{-1}\cdotp e
=\sum_{i=1}^n a_i f_i +y$ for some $a_i\in \C$ and $y \in \mfb. $
From the equality $$wb_2\cdotp e=\sum_{i=1}^n a_i f_i +y,$$ we
have $w\cdotp \Pi \subset (-\Pi) \cup \triangle^{+}.$ It follows
from Lemma $\ref{lem:1}$ that we have $w=w_J$ for some $J\subseteq
I.$ \ep
\begin{definition}
Fix  a dominant weight $\kappa$ of $\mfg,$ let $u\in(U^{+})^e $,
and define a regular function $\triangle_\kappa$ on $(U^+)^e$ as
$$\triangle_\kappa(u)=<u\cdotp v_{w_0\kappa},v_\kappa>,$$
where $v_\kappa$ and $v_{w_0\kappa}$ are highest and lowest weight
vectors in the representation $V_\kappa,$ respectively, and
$<u\cdotp v_{w_0\kappa},v_\kappa>$ is the coefficient of
$v_\kappa$ in the expansion of $u\cdotp v_{{w_0}\kappa}.$
\end{definition}
Note that the function $\triangle_\kappa$ is only defined up to a
choice of highest and lowest weight vectors.  The following two
lemmas and corollary are the generalizations of ideas in Lemma
$2.3$ of ~\cite{Riet1}.

\smallskip
\begin{lemma}~\label{lem:3}
Let  $u \in (U^{+})^e\cap B^{-}w^{J}B^{-}$ for some $J\subset I$,
and $\kappa$ be a dominant weight of $\mathfrak{g},$ then
$\triangle_\kappa(u)\neq 0$ precisely when \hs $w_0w_Jw_0\cdotp
v_{\kappa}=\pm v_{\kappa}.$
\end{lemma}
\begin{proof}
Write $u=b_{1}w_{0}w_{J}b_{2}$ for some $b_{1}\in B^{-}$ and
$b_{2}\in U^{-}$. Then we have $$u\cdotp
v_{w_0\kappa}=(b_{1}w_{0}w_{J}b_{2})\cdotp v_{w_{0}\kappa}= \pm(
b_{1}w_{0}w_{J})\cdotp v_{w_{0}\kappa}=\pm b_1w_0w_Jw_0\cdotp
v_{\kappa},$$ and hence  $<u\cdot v_{w_0 \kappa},v_{\kappa}>\ne0$
precisely if $w_{0}w_{J}w_{0}\cdotp v_{\kappa}=\pm v_{\kappa}.$
Indeed, suppose $w_{0}w_{J}w_{0}\cdot v_{\kappa}=\pm v_l$  for
some weight $l.$ Then unless $l$ is $\kappa,$ the action of
$B^{-}$ on $ v_l$ cannot have highest weight vector $v_\kappa$,
i.e., $<b\cdotp v_l,v_\kappa>=0$ for all $b\in B^{-},$ and hence
$<b_1\cdotp v_l,v_\kappa>=0.$
\end{proof}

Define an action of $w_0$ on $I$ by ${w_0}{i}=j$ if
${w_{0}}{s_i}{w_{0}}=s_j$. This action is well-defined since if
$w$ is an element of Weyl group $W$ and $s_\alpha$ is the
reflection in a root $\alpha$, then $w s_\alpha w
^{-1}=s_{w(\alpha)}$. Now $w_0$ is its own inverse, and $w_0$
sends a system of simple roots $\{\alpha_1,\dots,\alpha_n\}$ to
$\{-\alpha_1,\dots,-\alpha_n\}$, so that $w_0s_iw_0$ is $s_j,$
where $w_0(\alpha_i)=-\alpha_j.$

\smallskip
\begin{lemma}~\label{lem:4}
Fix $u\in (U^+)^e,$ and let $r_1,...,r_n$ be positive integers.
Define $J_u:=\{i\in I \mid \triangle_{{r_i}{\kappa_i}}(u)=0\}.$
Then $u$ is an element of $ B^{-}w^JB^{-}$ precisely when
$J=w_0J_u.$
\end{lemma}

\begin{proof}
Suppose $u\in B^{-}w^{J}B^{-}$. We will show that $ J=w_{0}J_{u}$.
By Lemma \ref{lem:3}, $i\notin J_{u}$, i.e.
$\triangle_{r_i{\kappa_i}}(u)\neq 0$ precisely if
$w_{0}w_{J}w_{0}\cdot v_{{r_i}{\kappa_i}}=\pm
v_{{r_i}{\kappa_i}},$ which is equivalent to saying
$$w_{0}w_{J}w_{0}\in \{w\in W\mid
w\cdotp({r_i}{\kappa_i})=r_i{\kappa_i}\}=\{w\in W|w\cdotp
\kappa_i=\kappa_i\}=<s_j|j\neq i>.$$ This is possible only when
$i\notin {w_0}J$. Therefore $J_u={w_0}J,$ i.e., $J=w_0J_u.$ For
the converse, suppose $ J=w_{0}J_{u}.$ Since $u$ is an element of
$(U^+)^e,$ by Lemma \ref{lem:2}, $u\in B^{-}w^{J^{\prime}}B^{-} $
for some $J^{\prime}\subset I.$ Applying the same argument as
above, we get $J^\prime=w_0 J_u.$ Since $J=w_0J_u,$ we get
$J^\prime=J.$
\end{proof}
 Note that from the proof the set $J_u$ does not depend on the choices of $r_1,...,r_n.$

\smallskip
\begin{corollary}~\label{cor:1}
As a set, $\overline{B^{-}w^{J}B^{-}}\cap (U^+)^e=\{u\in (U^+)^e|
\triangle_{{r_i}{\kappa_i}}(u)=0 $ for all $i\in w_0J$ \}.
\end{corollary}
\begin{proof}
 $u\in \overline{B^{-}w^{J}B^{-}}\cap (U^+)^e $ precisely if
 $u\in (B^{-}w^{J^{\prime}}B^{-})\cap (U^{+})^e$ for
some $J^{\prime}$ with $J\subseteq J^{\prime}\subseteq I$ if and
only if  $w_0J^{\prime}=J_u$  for some $J^{\prime}$ with
$J\subseteq J^{\prime}\subseteq I$ by Lemma \ref{lem:4}.  By the
definition of $J_u,$ this is equivalent to saying that
$\triangle_{{r_i}{\kappa_i}}(u)=0 $ for all $i\in w_0J.$
\end{proof}

In the rest of the section, using Corollary \ref{cor:1}, we
compare two presentation of the quantum cohomology rings of
$\ogo,$ $\lg$ and $\oge.$ Let \be\label{M(m)} M(m):=
\left(\begin{array}{cccccc}
1&X_1&X_2&\cdots &X_{m-2}&X_{m-1}\\
&1&X_1& & &X_{m-2}\\
& &\dd&\dd& &\vd\\
&&&\dd&X_1&X_2\\
&&&&1&X_1\\
&&&&&1
\end{array}\right),
\ee and let \be \label{2 paffian}
X_{i,i}:=X_i^2+2\sum_{k=1}^i(-1)^k X_{i+k}X_{i-k},\ee
 where
$X_0=1.$
\subsection{Case of $\ogo$}
We fix the nilpotent element $e$ of $\mathfrak{sp}_{2n}(\C)$;
$$e=\sum_{i=1}^n e_{i}=\sum_{i=1}^{2n-1}E_{i,i+1}.$$

The elements $u$ of $(U^+)^e$ in $Sp_{2n}(\C)$ are exactly
matrices of the form $M(2n)$ in (\ref{M(m)}) together with the
relations, $X_{i,i}=0,$ $i=1,...,n-1.$ Therefore
 we can identify the coordinate ring of $(U^+)^e$
with a quotient of $\C[X_1,...,X_{2n-1}]$ modulo the relations
$X_{i,i}=0, i=1,...,n-1.$

\begin{definition}
We fix the subset $J=\{1,...,n-1\}\subset \{1,2,...,n\}=I$. Define
$\mcv_n^C$ to be a closed subvariety $$\mcv_{n}^C:=\U+e \cap
\overline{B^{-}w^{J}B^{-}}\subset \U+e \subset Sp_{2n}(\C).$$
\end{definition}
We note that $\mcv_{n}^C$ is a $1$-dimensional closed subvariety
of $(U^+)^e$ because it is cut out by $(n-1)$ equations by
Corollary \ref{cor:1} and the dimension of $(U^+)^e$ is $n$ .

\begin{theorem}[Peterson]\label{thm;ogo} Let $\mathfrak{I}_o$ be
the ideal of the ring $\C[\tau_1,...,\tau_n]$ generated by
$\tau_{i,i}$ for $i=1,...,n-1$ , and let
$\mathfrak{R}_o:=\C[\tau_1,...,\tau_n]/\mathfrak{I}_o.$ Then the
map $\phi_C:\mathfrak{R}_o \rightarrow \mo(\mcv_{n}^C)$ defined by
$\tau_i \mapsto \frac{1}{2} X_i$ for $i=1,...,n$ is an
isomorphism.
\end{theorem}

\begin{remark}
Note that the the ring $\mathfrak{R}_o$ is isomorphic to quantum
cohomolgy ring $qH^*_{\C}(\oge)$ given in \ref{subsec:ogo} which
is in turn isomorphic to $qH^*_{\C}(\ogo)$. Therefore Theorem
\ref{thm;ogo} implies that $qH^*_{\C}(\ogo)$ is isomorphic to the
reduced coordinate ring $\mo(\mcvc)$. The next lemma is a part of
the proof of this theorem, and, in $\ref{completion of proofs},$
we will complete the proof of the theorem by showing that
$\mathfrak{R}_o$ is reduced.
\end{remark}

\begin{lemma}\label{lemma new1}
The map from $\C[\tau_1,...,\tau_n]$ to $\mo(\mcvc)$ defined by
$\tau_i \mapsto \frac{1}{2}X_i $ for  $i=1,...,n$  is a surjective
ring homomorphism with the kernel the radical of the ideal
generated by $\tau_{i,i},$ $i=1,...,n-1,$ with $\tau_{i,i}$ given
in $\ref{subsec:ogo}.$
\end{lemma}

\begin{proof}
Since $J=w_0\cdotp J=\{1,...,n-1\}$, by Corollary $\ref{cor:1},$
$$\mcvc=\{u\in \U+e \mid \triangle_{\kappa_r}(u)=0,
r=1,2,...,n-1\}.$$ For each
 $r=1,...,n-1$, $$\kappa_r=l_1+\cdots + l_r  \hspace{0.08in}
\textrm{and} \hspace{0.08in} w_0\cdotp \kappa_r=-l_1-\cdots
-l_r,$$ and their weight vectors are $$\xi_{1}\wedge\cdots\wedge
\xi_{r} \hspace{0.08in} \textrm{and} \hspace{0.08in}
\xi_{2n+1-r}\wedge\cdots\wedge \xi_{2n},$$ respectively. Therefore
$\triangle_{\kappa_r}(u)$ is the determinant of $(r\times r)$
submatrix on the right upper corner of the matrix $M(2n).$  The
vanishing of the $\triangle_{\kappa_r}(u)$, $r=1,...,n-1$,
inductively implies that the coordinates $X_{2n-1},...,X_{n+1}$
vanish. This proves the lemma.
\end{proof}

\subsection{Case of $\lg$}
We recall that the special result of Peterson holds for minuscule
parabolic subgroups, and $P_n \subset Sp_{2n}(\C)$ is not
minuscule (but cominuscule). But we can still find an analogue of
Peterson's result for this case.
 For simplicity, we take $e^{\prime}$ of $\mathfrak{so}_{2n+1}(\C)$, which
 is a normalization of the principal nilpotent element $e$ given
 in \ref{odd orthogonal group},
$$e^{\prime}:=\sum_{i=1}^{n}(E_{i,i+1}+E_{i^{*}-1,i^{*}})=\sum_{i=1}^{2n}E_{i,i+1},$$
where $i^{*}=2n+2-i.$ We note that all the lemmas in \ref{Bruhat
cells 1} hold for $e^{\prime}.$ The elements $u$ of
$(U^+)^{e^{\prime}}$ in $SO_{2n+1}(\C)$ are precisely matrices of
the form $M(2n+1)$ in (\ref{M(m)}) together with $ X_{i,i}=0,$
$i=1,...,n.$

Therefore the coordinate ring $(U^+)^{e^{\prime}}$ can be
identified with a quotient of the polynomial ring
$\C[X_1,...,X_{2n}]$ modulo the relations $X_{i,i}=0,$
$i=1,...,n.$

\begin{definition}
We fix the subset $J=\{1,...,n-1\}\subset \{1,2,...,n\}=I$. Define
$\mcv_n^B$ to be a closed subvariety
$$\mcv_n^B:=(U^+)^{e^{\prime}}\cap\overline{B^{-}w^{J}B^{-}} \subset (U^+)^{e^{\prime}} \subset
SO_{2n+1}(\C).$$
\end{definition}

\begin{theorem}\label{2nd}
Let $\mathfrak{I}_L$ be the ideal of the polynomial ring
$\C[\Qt_1,...,\Qt_{n+1}]$ generated by $\Qt_{i,i}$ for
$i=1,...,n,$ and let
$\mathfrak{R}_L:=\C[\Qt_1,...,\Qt_{n+1}]/\mathfrak{I}_L$. Then the
map $\phi_B:\mathfrak{R}_L \rightarrow \mo(\mcvb)$ defined by
$\Qt_i \mapsto X_i$ for $i=1,...,n+1$ is an isomorphism. Here
$\Qt_i$ are understood as $\Qt_i(X_{n+1}).$
\end{theorem}

\begin{remark}
  We note that the map
from $qH^*_\C(\lg)$ to $\mathfrak{R}_L$ defined by
\begin{displaymath} \left.\begin{array}{cccc}
\sigma_i&\mapsto& \Qt_i & \textrm{for} \hs i=1,...,n,\\
q&\mapsto& 2\Qt_{n+1}
\end{array}\right.
\end{displaymath}
is an isomorphism (\cite{KT2}), and so Theorem \ref{2nd} implies
that the quantum cohomology ring $qH^*_\C(\lg)$ is isomorphic to
the reduced coordinate ring $\mo(\mcvb).$ The following lemma and
\ref{completion of proofs} give a full proof of this theorem.
\end{remark}

\begin{lemma}\label{lemma new2}
The map from $\C[\Qt_1,...,\Qt_{n+1}]$ to $\mo(\mcvb)$ defined by
$\Qt_i \mapsto X_i$ for $i=1,...,n+1$ is a surjective ring
homomorphism with the kernel the radical of the ideal generated by
$\Qt_{i,i},$ $i=1,...,n.$
\end{lemma}

\bp In this case, $J=w_0\cdotp J=\{1,...,n-1\},$
$$\triangle_{\kappa_r}(u)=<u\cdotp \xi_{2n+2-r}\wedge\cdots \wedge
\xi_{2n+1},\xi_1\wedge\cdots \wedge \xi_r>,$$ and $$\mcvb=\{u\in
(U^+)^{e^{\prime}}\mid \triangle_{\kappa_r}(u)=0, r=1,...,n-1\}.$$
Therefore, as in the case of $\ogo,$ $X_{2n},...,X_{n+2}$
inductively vanish. So the lemma follows. \ep

\subsection{Case of $\oge$}
In this case we take the principal nilpotent element of
$\mathfrak{so}_{2n+2}(\C)$
$$e=\sum_{i=1}^{n+1}e_i=\sum_{i=1}^{n}(E_{i,i+1}-E_{2n+2-i,2n+3-i})+(E_{n,n+2}-E_{n+1,n+3}).$$
The elements $v$ of $(U^+)^e$ in $SO_{2n+2}(\C)$ are exactly
matrices of the
 block form \be\label{soe:u} v=\left(\begin{array}{cc}
A&B\\
O&C
\end{array}\right), \ee satisfying the following relations $(\ref{7.2})$ and
$(\ref{7.3})$,
\begin{equation}\label{7.2}
X_{r,r}^{\prime}:={X_r^{\prime}}^2-2\sum_{i=1}^r(-1)^{i+1}{X_{r+i}^{\prime}}{X_{r-i}^{\prime}}=0
\trm r=1,...,n-1,
\end{equation}
where $X_k^{\prime}=X_k \hspace{0.08in}\textrm{if}
\hspace{0.08in}k\leq n-1,$ and $X_k^{\prime}=Y_k
\hspace{0.08in}\textrm{ otherwise},$
\begin{equation}\label{7.3}
X_n^2-2\sum_{i=0}^nY_{n+i}X_{n-i}=0.\end{equation} The submatrices
of $v$ are given by
\begin{displaymath} A=\left(\begin{array}{cccccc}
1&X_1&X_2&\cdots&X_{n-1}&X_n\\
0&1&X_1&X_2&\cdots&X_{n-1}\\
0&0&1&X_1&\dd&\vd\\
\vd&\vd&0&\dd&\dd&X_2\\
0&0&0&\dd&1&X_1\\
0&0&0&\cdots&0&1
\end{array}\right),
\end{displaymath}

\begin{displaymath}B= \left(\begin{array}{c|ccccc}
2Y_n-X_n&-2Y_{n+1}&2Y_{n+2}&\cdots&(-1)^{n-1}2Y_{2n-1}&(-1)^n2Y_{2n}\\
X_{n-1}&-2Y_{n}&2Y_{n+1}&\cdots&(-1)^{n-1}2Y_{2n-2}&(-1)^n2Y_{2n-1}\\
X_{n-2}&-2X_{n-1}&2Y_n&\cdots&(-1)^{n-1}2Y_{2n-3}&(-1)^n2Y_{2n-2}\\
X_{n-3}&-2X_{n-2}&2X_{n-1}&\cdots&(-1)^{n-1}2Y_{2n-4}&(-1)^n2Y_{2n-3}\\
\vdots&\vdots&&\vdots&\vdots&\vdots\\
X_3&-2X_4&2X_5&\cdots&(-1)^{n-1}2Y_{n+2}&(-1)^{n}2Y_{n+3}\\
X_2&-2X_3&2X_4&\cdots&(-1)^{n-1}2Y_{n+1}&(-1)^n2Y_{n+2}\\
X_1&-2X_2&2X_3&\cdots&(-1)^{n-1}2Y_n&(-1)^n2Y_{n+1}\\
\hline
 0&-X_1&X_2&\cdots&(-1)^{n-1}X_{n-1}&W_n

\end{array}\right),
\end{displaymath}

\begin{displaymath}C= \left(\begin{array}{ccccccccccc}
1&-X_1&X_2&\cdots&(-1)^{n-1}X_{n-1}&Z_n\\
0&1&-X_1&X_2&\dd&(-1)^{n-1}X_{n-1}\\
0&0&1&-X_1&\dd&\vd\\
\vd&0&0&\dd&\dd&X_2\\
0&\vd&\dd&\dd&1&-X_1\\
0&0&\cdots&0&0&1
\end{array}\right),
\end{displaymath}

and  the submatrix $O$ is the $(n+1)\times (n+1)$-zero matrix.
Here
\begin{displaymath}W_n= \left\{\begin{array}{cc}
X_n-2Y_n & \textrm{if n is odd,}\\
X_n &  \textrm{if n is even.}
\end{array}\right.
\end{displaymath}
and
\begin{displaymath}Z_n= \left\{\begin{array}{cc}
-X_n &  \textrm{if n is odd,}\\
2Y_n-X_n &  \textrm{if n is even.}
\end{array}\right.
\end{displaymath}

The coordinate ring of  $\U+e$ can be identified with a quotient
of the ring  $\C[X_1,...,X_n,Y_n,...,Y_{2n}]$ modulo the relations
$(\ref{7.2})$ and $(\ref{7.3})$.
\begin{definition}
For the subset $J=\{1,...,n\}\subset \{1,...,n+1\}=I,$ we define
$\mcv_{n+1}^D$ to be a closed subvariety
$\mcv_{n+1}^D:=\U+e\cap\overline{B^{-}w^{J}B^{-}}\subset \U+e
\subset SO_{2n+2}(\C).$
\end{definition}

\begin{theorem}[Peterson]\label{thm:oge}
There is an isomorphism $\phi_D:\mathfrak{R}_o \rightarrow
\mo(\mcvd)$ that takes $\tau_i$ to $X_i^{\prime},$ $i=1,...,n.$
\end{theorem}
Note that since the ring $\mathfrak{R}_o$ is isomorphic to the
quantum cohomology ring $qH_\C^*(OG^e(n)),$ this theorem implies
that $qH_\C^*(OG^e(n))$ is isomorphic to the reduced coordinate
ring $\mo(\mcvd).$ The following lemma is a part of the proof of
the theorem. To complete the proof of the theorem, it remains to
show that $\mathfrak{R}_o$ is reduced, which will be done in
$\ref{completion of proofs}.$
\begin{lemma}\label{lemma new3}
The map from $\C[\tau_1,...,\tau_n]$ to $\mo(\mcvd)$ defined by
$\tau_i \mapsto X_i^{\prime}$ for $i=1,...,n$ is a surjective ring
homomorphism with the kernel the radical of the ideal generated by
$\tau_{i,i},$ $i=1,...,n-1.$
\end{lemma}
\bp First recall from \ref{sec:basic D} that the fundamental
weights of $SO_{2n+2}(\C)$ consist of $(n+2)$-weights, $\omega_r$
for $r=1,...,n,$ $\omega_{n+1}^{-}=2\kappa_{n}$ and
$\omega_{n+1}^{+}=2\kappa_{n+1}.$

For the sake of simplicity, let \begin{displaymath}\omega_{n+1}:=
\left\{\begin{array}{ccc}
\omega_{n+1}^{-}& \textrm{if}\hs n  \hs \textrm{is} \hs\textrm{odd},\\
\omega_{n+1}^{+} &\textrm{if}\hs n \hs \textrm{is}\hs\textrm{
even}.
\end{array}\right.
\end{displaymath}
 By (\ref{w0}), for $J=\{1,...,n\},$ we have
\begin{displaymath}w_0J= \left\{\begin{array}{ccc}
\{1,...,n\}& \textrm{if}\hs n  \hs \textrm{is} \hs\textrm{odd},\\
\{1,...,n-1,n+1\} &\textrm{if}\hs n \hs \textrm{is}\hs\textrm{
even}.
\end{array}\right.
\end{displaymath}
Therefore Corollary \ref{cor:1} implies that, for any $n$, odd or
even, the variety $\mcv_{n+1}^D$ is defined by the equations
$\triangle_{\omega_{k}}(v)$ for $k=1,...,n-1,n+1,$ and the
coordinate ring $\mo(\mcv_{n+1}^D)$ is isomorphic to the ring
\begin{equation}\label{7.5}
\C[X_1,..,X_n,Y_n,...,Y_{2n}] /\textrm{Rad}(\ma),
\end{equation}
 where $\textrm{Rad}(\ma)$ is the radical of  ideal $\ma$, and $\ma$ is generated by the polynomials
 in $(\ref{7.2})$ and $(\ref{7.3})$, and polynomials $\triangle_{\omega_k}(v)$ for $k\in\{1,...,n-1,n+1\}.$
From the relation (\ref{b.1}) in $\ref{sec:basic D}$ and the fact
that $ \triangle_{\omega_k}(v)$ lie in
 $\textrm{Rad}(\ma)$ for $k\in\{1,...,n-1,n+1\},$
the polynomial  $\triangle_{\omega_n}^2(v)$ lies in $\textrm{
Rad}(\ma)$, and so does $\triangle_{\omega_n}(v).$  For each
$k=1,...,n,$ the matrix coefficient $\triangle_{\omega_k}(v)$ is
given by the determinant of the $(k\times k)$ submatrix of $B$ on
the right upper corner. Note that if $n$ is odd, then the weights
$$\omega_{n+1}=l_1+\cdots +l_n-l_{n+1} \hspace{0.08in}
\textrm{and} \hspace{0.08in
}w_0\omega_{n+1}=-l_1-\cdots-l_n+l_{n+1}$$  have weight vectors
$$\xi_1\wedge \cdots \xi_n\wedge \xi_{n+2} \hspace{0.08in}
\textrm{and} \hspace{0.08in} \xi_{n+1}\wedge \xi_{n+3}\wedge
\cdots \wedge \xi_{2n+2},$$ respectively, and if $n$ is even, the
weights
$$\omega_{n+1}=l_1+\cdots +l_n+l_{n+1}  \hspace{0.08in}
\textrm{and} \hspace{0.08in}
w_0\omega_{n+1}=-l_1-\cdots-l_n+l_{n+1}$$ have weight vectors
$$\xi_1\wedge \cdots \wedge \xi_n\wedge \xi_{n+1} \hspace{0.08in}
\textrm{and} \hspace{0.08in} \xi_{n+1}\wedge \xi_{n+3}\wedge
\cdots \wedge \xi_{2n+2},$$ respectively. Therefore if $n$ is odd,
$\triangle_{\omega_{n+1}}(v)$ is the determinant of a submatrix
$M_n$ defined below of $v$ with $(1,...,n,n+2)$-th rows and
$(n+1,n+3,...,2n+2)$-th columns, and if $n$ is even,
$\triangle_{\omega_{n+1}}(v)$ is the determinant of the submatrix
$M_n$ of $v$ with  $(1,...,n,n+1)$-th rows and
$(n+1,n+3,...,2n+2)$-th columns. Here the submatrix $M_n$ is given
as follows.
\begin{displaymath}
M_n= \left(\begin{array}{c|cccccc}
X_n&-2Y_{n+1}&2Y_{n+2}&\cdots& (-1)^{n-1}2Y_{2n-1}&(-1)^n2Y_{2n}\\
X_{n-1}&-2Y_{n}&2Y_{n+1}&\cdots& (-1)^{n-1}2Y_{2n-2}& (-1)^n2Y_{2n-1}\\
X_{n-2}&-2X_{n-1}&2Y_n&\cdots& (-1)^{n-1}2Y_{2n-3}& (-1)^n2Y_{2n-2}\\
X_{n-3}&-2X_{n-2}&2X_{n-1}&\cdots& (-1)^{n-1}2Y_{2n-4}& (-1)^n2Y_{2n-3}\\
\vdots & \vdots&\vdots&&\vdots&\vdots\\
X_3&-2X_4&2X_5 &\cdots&(-1)^{n-1}2Y_{n+2}& (-1)^{n}2Y_{n+3}\\
X_2&-2X_3&2X_4&\cdots& (-1)^{n-1}2Y_{n+1}& (-1)^n2Y_{n+2}\\
X_1&-2X_2&2X_3&\cdots& (-1)^{n-1}2Y_n& (-1)^n2Y_{n+1}\\
\hline
 a&-X_1&X_2&\cdots& (-1)^{n-1}X_{n-1}& P_n
\end{array}\right),
\end{displaymath}
where
\begin{displaymath}P_n= \left\{\begin{array}{cc}
Z_n& \textrm{if n is odd,}\\
W_n &  \textrm{if n is even,}
\end{array}\right.
\end{displaymath}
and
\begin{displaymath}a= \left\{\begin{array}{cc}
0 &  \textrm{if n is odd,}\\
1 &  \textrm{if n is even.}
\end{array}\right.
\end{displaymath}

The containment of $\triangle_{\omega_k}(v)$ in
$\textrm{Rad}(\ma)$ for $k=1,...,n$ implies inductively that
$Y_{2n},...,Y_{n+1}$ lie in the ideal $\textrm{Rad}(\ma)$.
Therefore the ring in (\ref{7.5}) is isomorphic to
\begin{equation}\label{7.6}
\C[X_1,...,X_n,Y_n]/\textrm{Rad}(\mb).
\end{equation}
Here the ideal $\mb$ is generated by the following polynomials
$(\ref{p0}),$  $(\ref{p2})$ and $(\ref{p3})$;

\begin{equation}\label{p0}
{X_r^{\prime}}^2-2\sum_{i=1}^{r}(-1)^{i+1}{X_{r+i}^{\prime}}{X_{r-i}^{\prime}}
\trm r=1,...,n-1,
\end{equation}

\be\label{p2} X_n^2-2X_nY_n, \ee

\be \label{p3}
\triangle_{\omega_{n+1}}(v)\av_{Y_{2n}=\cdots=Y_{n+1}=0}.\ee

From the matrix $M_n,$ we easily see that
$\triangle_{\omega_{n+1}}(v)\av_{Y_{2n}=\cdots=Y_{n+1}=0}=(-1)^n2^{n-1}X_n^2Y_n^{n-1}$
for  $n,$  both odd and even. From  $(\ref{p2})$ and $(\ref{p3}),$
it follows that in the ring
$\C[X_1,...,X_n,Y_n]/\textrm{Rad}(\mb),$ we have the relations
$$X_n^2=2X_n Y_n \hs \textrm{and}\hs X_n^2 Y_n^{n-1}=0.$$
This implies
$$X_n^{2n+1}=2^nX_n^{n+1}Y_n^n=2^n(X_n^{n-1}Y_n)(X_n^2Y_n^{n-1})=0.$$
Since the ring $\C[X_1,...,X_n,Y_n]/\textrm{Rad}(\mb)$ is reduced,
we have $X_n=0$  in  the ring
  $\C[X_1,...,X_n,Y_n]/\textrm{Rad}(\mb).$ Therefore the coordinate
ring $\mo(\mcvd)$ is isomorphic to the following reduced ring \be
\label{7.7} \C[X_1^{\prime},...,X_{n}^{\prime}]/\textrm{Rad}(\mc),
\ee where the ideal $\mc$ is generated by the polynomials
\begin{equation}
{X_r^{\prime}}^2-2\sum_{i=1}^{r}(-1)^{i+1}{X_{r+i}^{\prime}}{X_{r-i}^{\prime}}
\trm r=1,...,n-1.
\end{equation}
This proves the lemma. \ep

\begin{remark}
From Lemma \ref{lemma new1}, we have a well-defined map (a priori
not an isomorphism ) $\phi_C: \mathfrak{R}_o\rightarrow
\mathcal{O}(\mcvc).$ So each element $F$ of $\mathfrak{R}_o$
induces a function $\phi_C(F)$ on $\mcvc.$ We will write $\dot{F}$
for $\phi_C(F)$. Similarly, $\phi_B(G):=\dot{G}$ and
$\phi_D(H)=\dot{H}$ can be considered as functions on $\mcvb$ and
$\mcvd$ for any $G \in \mathfrak{R}_L$ and $H\in \mathfrak{R}_o.$
\end{remark}

\section{ELEMENTS OF VARIETIES $\mcvc,$ $\mcvb$ AND $\mcvd$}
\label{sec:5} In this section, we describe elements of $\mcvc,$
$\mcvb$ and $\mcvd.$ To do this, we need the following definitions
and notations.
\subsection{Definitions and Notations.}\label{def and nota}
Let us call $I=(i_1,...,i_n)\in \mathcal{T}_n$ $exclusive$  if
$\zeta^{i_k}\ne -\zeta^{i_l}$ for all $k,l=1,...,n.$ So if $I$ is
exclusive, then the set $\{\zeta^{i_1},..., \zeta^{i_n}\}$ can not
have both of two antipodal points, but only one.  We also call
$I=(i_1,...,i_n)\in \mathcal{T}_n$ $self$-$symmetric$ if the set
$\{\zeta^{i_1},..., \zeta^{i_n}\}$ contains both $\zeta^{i_k}$ and
the complex conjugate $\overline{\zeta^{i_k}}.$  Recall that
$\mathcal{I}_n$ denotes the set of all exclusive $n$-tuples
$I=(i_1,...,i_n).$ Denote by $\In^s$ the set of exclusive and
self-symmetric $n$-tuples $I=(i_1,...,i_n)$ in $\mathcal{T}_n.$ We
can easily check that $|\mathcal{I}_n|=2^n=|\mathcal{D}(n)|,$ and
$|\In^s|=2^{\lfloor\frac{n}{2}\rfloor}.$

For $a_1,...,a_n \in \C,$ define a matrix
$\tilde{u}_m(a_1,...,a_n)\in SL_{m}(\C)$ as \be \label{special
matrix}
 \tilde{u}_m(a_1,...,a_n):= \left(\begin{array}{cccccccc}
1&a_1&a_2&\cdots&a_n&0&\cdots&0\\
&1&a_1&\dd&&\dd&&\vdots\\
&&1&\dd&&&\dd&0\\
&&&\dd&&&&a_n\\
&&&&\dd&\dd&\dd&\vd\\
&&&&&1&a_1&a_2\\
 &&&&&&1&a_1 \\
 &&&&&&&1\\
\end{array}\right), \ee

For $x_1,...,x_n \in \C,$ let
$u_{m}(x_1,...,x_n):=\tilde{u}_m(E_1(x_1,...,x_n),...,E_n(x_1,...,x_n)).$
Similarly, let\\
$\tilde{v}_{2n+2}(a_1,...,a_n)$ be the matrix $v$ in (\ref{soe:u})
with $X_i^{\prime}$ replaced by $a_i$ for $i=1,...,n,$ and with
$X_n$ and other $X_i^{\prime}$ replaced by $0.$ Let
$v_{2n+2}(x_1,...,x_n):=\tilde{v}_{2n+2}(E_1(x_1,...,x_n),...E_n(x_1,...,x_n)).$

We suppress the subscripts in $u_m$ and $v_{2n+2}$ if any
confusion does not arise.
\subsection{Elements of varieties $\mcvc,$  $\mcvb$ and $\mcvd$.}
We now characterize the elements of $\mcvc,$ $\mcvb$ and $\mcvd$
in terms of $x_i.$
\begin{lemma}\label{lemma:characterization}
For $I=(i_1,...,i_n)\in \mathcal{T}_n,$ $I$ is exclusive precisely
when $E_{i}(\zeta^{2I})=0$ for all $i=1,...,n-1.$
\end{lemma}

\begin{proof}
Obviously,
$I_0:=(-\frac{n-1}{2},-\frac{n-1}{2}+1,...,\frac{n-1}{2})$ is
exclusive, $E_i(\zeta^{2I_0})=0$ for $i=1,...,n-1$ and
$E_n(\zeta^{2I_0})=1.$ If $I$ is exclusive, then $2I=2I_0$ up to
order, and so $E_i(\zeta^{2I})=0$ for $i=1,...,n-1.$ For the
converse, suppose
 $E_{i}(\zeta^{2I})=0$ for $i=1,...,n-1.$ If $I$ is not exclusive,
 then  there are $k$ and $l$ such that $\zeta^{2i_{k}}=\zeta^{2i_l}.$
Consider the following polynomial in $z$, of which zeros are
$\eta_m:=\zeta^{2i_m}$ for $m=1,...,n,$
$$g(z)=\prod_{i=1}^n(z-\eta_i).$$
Since $E_i(\eta_1,...,\eta_n)=0$ for $i=1,...,n-1,$ we have

\be \label{eee}
\prod_{i=1}^n(z-\eta_i)=z^n+(-1)^nE_n(\eta_1,...,\eta_n).\ee If we
take derivatives and evaluation at $z=\eta_k$ of both sides of
$(\ref{eee})$, the left hand side is zero,  and the right hand
side is $n\eta_k^{n-1},$ which is nonzero, which yields a
contradiction.
\end{proof}

\begin{lemma}\label{pp}
The matrix $u_{2n}(x_1,...,x_n)$ is an element of $\mcvc$ if and
only if $(x_1,...,x_n)$ can be written (up to order) as
$(x_1,...,x_n)=t\zeta^I$ for some $I \in \mathcal{I}_n$ and $t\in
\C.$ All elements of $\mcvc$ are of the form $u_{2n}(t\zeta^I)$
with
 $t\in \C$ and $I \in \In.$
\end{lemma}

\begin{proof}
First note that from Lemma \ref{lemma new1} all the elements of
$\mcvc$ are of the form $u_{2n}(x_1,...,x_n)$ for complex numbers
$x_1,...,x_n$, and, by the definitions of $X_{i,i}$ and
$u_{2n}(x_1,...,x_n),$ we have
$$X_{i,i}(u_{2n}(x_1,...,x_n))=\Qt_{i,i}(x_1,...,x_n),$$  which is equal to
$E_{i}(x_1^2,...,x_n^2)$ by the property $(1)$ of
$\Qt$-polynomials in \ref{subsec:sym poly}.
 Therefore  $\mcvc$  consists of all $u_{2n}(x_1,...,x_n)$
such that $E_{i}(x_1^2,...,x_n^2)=0,$ $i=1,...,n-1.$ To prove the
lemma, it suffices to show that $E_{i}(x_1^2,...,x_n^2)=0$ for
$i=1,...,n-1$ if and only if $(x_1,...,x_n)=t\zeta^I$ for some $I
\in \mathcal{I}_n$ and $t\in \C.$ This hold trivially for
$(x_1,...x_n)=(0,...,0).$ So we assume $(x_1,...x_n)\ne(0,...,0).$
First suppose that $E_i(x_1^2,...,x_n^2)=0$ for all $i=1,...,n-1$.
Consider the following polynomial in $z,$
$$ f(z)=\prod_{i=1}^n(z-x_i^2).$$
Since $E_i(x_1^2,...,x_n^2)=0$ for $i=1,...,n-1$, we have
\begin{equation}\label{ele}
\prod_{i=1}^n(z-x_i^2)=z^n+{(-1)^n}E_n(x_1^2,..,x_n^2).
\end{equation}
But since $x_1^2,...,x_n^2$ are roots of $f(z)$,  we have \be
\label{aa} x_1^{2n}=\cdots
=x_n^{2n}=(-1)^{n+1}E_n(x_1^2,...,x_n^2).\ee

 If
$x_i=0$ for some $1\leq i\leq n,$ then trivially all $x_i$ are
zero $(\ref{aa})$, which contradict the assumption that
$(x_1,...x_n)\ne(0,...,0),$ and hence  all $x_i$ are nonzero.
Furthermore  all $x_i$ are distinct. Indeed if all $x_i$ are not
distinct, say, if $x_1=x_2$, taking derivatives and evaluation at
$z=x_1^2$ of both sides of $(\ref{ele})$, the left hand side is
zero , and the right hand side is $nx_1^{2n-2},$ which is nonzero,
and hence we get a contradiction.  Now by the description
$(\ref{aa})$ of $x_1,...,x_n$, we can write $(x_1,...,x_n)=t
\zeta^I$ for some $I\in \mathcal{T}_n$ and nonzero $t\in \C.$ But
since
$$E_i(x_1^2,...,x_n^2)=E_i(t^2\zeta^{2I})=t^{2i}E_{i}(\zeta^{2I})=0,$$
and $t$ is nonzero,
 we have  $E_{i}(\zeta^{2I})=0,$  and so, by Lemma
$\ref{lemma:characterization},$  $I$ is exclusive. The other
direction is trivial by Lemma $\ref{lemma:characterization}.$ This
completes the proof.
\end{proof}

\begin{lemma}\label{ppa}
The matrix $u_{2n+1}(x_1,...x_{n+1})$ is an element of $\mcvb$ if
and only if $(x_1,...,x_{n+1})$ can be written (up to order) as
$(x_1,...,x_{n+1})=t\zeta^I$ for some  $I\in \mathcal{I}_{n+1}$
and $t\in \C.$ All elements of $\mcvb$ are of this form
$u_{2n+1}(t\zeta^I)$ with $t\in \C$ and $I \in \mathcal{I}_{n+1}$.
\end{lemma}
\begin{proof}
 The same argument as in the proof of Lemma $\ref{pp}$ applies,
 with $``u_{2n}(x_1,...,x_n)"$ and $``X_{i,i}=0$
 for $i=1,...,n-1"$ replaced by $``u_{2n+1}(x_1,...,x_{n+1})"$ and $``X_{i,i}=0$
 for $i=1,...,n",$ respectively.

\end{proof}
\begin{lemma}\label{ppb}
The matrix $v_{2n+2}(x_1,...x_{n})$ is an element of $\mcvd$ if
and only if $(x_1,...,x_{n})$ can be written (up to order) as
$(x_1,...,x_{n})=t\zeta^I$ for some  $I\in \mathcal{I}_{n}$ and
$t\in \C.$ All elements of $\mcvd$ are of this form
$v_{2n+2}(t\zeta^I)$ with $t\in \C$ and $I \in \mathcal{I}_{n}$.
\end{lemma}
\begin{proof}
The argument of the proof of Lemma $\ref{pp}$ applies, with
$``u_{2n}(x_1,...,x_n)"$ and $``X_{i,i}=0$
 for $i=1,...,n-1"$ replaced by $``v_{2n+2}(x_1,...,x_n)"$ and $``X_{i,i}^{\prime}=0$
 for $i=1,...,n-1",$ respectively.
\end{proof}
The evaluations of  the functions $\dot{q}$ defined by the quantum
variables will be used to find the Vafa-Intriligator type formulas
in Section \ref{sec:7}.

\begin{lemma}\label{lemma:quant12}
\bn
\item For the quantum variable $q$ for $\ogo$ and
$u=u_{2n}(t\zeta^I)\in \mcvc,$ we have
$$\dot{q}(u)=\frac{1}{4}t^{2n}.$$
\item For the quantum variable $q$ for $\lg$ and
$u=u_{2n+1}(t\zeta^I)\in \mcvb,$ we have
$$\dot{q}(u)=2t^{n+1}E_{n+1}(\zeta^I), \hspace{0.06in}
\textrm{and}\hspace{0.1in} \dot{q}^2(u)=4t^{2n+2}.$$\en
\end{lemma}
\begin{proof}
$(1)$ follows directly from the fact that  $q=\tau_n^2$ correspond
to $\Pt_n^2$ by \ref{subsec:ogo}, and
$$\Pt_{n,n}(x_1,...,x_n)=\frac{1}{4}E_n(x_1^2,...,x_n^2)=(-1)^{n+1}\frac{1}{4}x_i^{2n}.$$
Since  $q$ corresponds to $2 E_{n+1},$ we get $(2)$
\end{proof}
\subsection{The varieties $\mcvb,$ $\mcvc$ and $\mcvd$}
Let $$U_{2n}:=\{z\in \C\mid z^{2n}=1\}=\{\zeta^i\mid
i=0,...,2n-1\}.$$ Consider the natural actions of the groups
$U_{2n}$ and $\Z_2$ on the sets $\{\zeta^I\mid I \in \In \}$ and
$\{\zeta^I\mid I \in \In^s \},$ respectively.  Let $O_n$ and $o_n$
be the sets of orbits for the actions $U_{2n}$ and $\Z_2,$
respectively. Then we count $$|O_n|=\frac{1}{2n}\sum_{d| n \atop d
\hspace{0.04in}\textrm{odd} }\varphi(d)2^{\frac{n}{d}},
\hspace{0.08in}\textrm{and}\hspace{0.08in}|o_n|=2^{\lfloor\frac{n}{2}\rfloor-1},$$
where $\varphi(d)$ is the number of positive integers $k\leq d$
such that $k$ is relatively prime to $d.$ For $I\in \In, $ let
$\C^{*}_I:=\{u_{2n}(t\zeta^I)\in \mcvc \mid t\in \C^{*}\},$ and
for $I \in \In^s,$ let $\R^{*}_I:=\{u_{2n}(t\zeta^I) \in \mcvc
\mid t\in \R^{*}\}.$ Then we have
$$\mcvc=\{id\} \bigcup_{I\in \In}\C^{*}_I, \hspace{0.08in} \textrm{and}\hspace{0.08in}
\mcvc(\R)=\{id\} \bigcup_{I\in \In^s}\R^{*}_I,$$ where $\mcvc(\R)$
is the set of real points in $\mcvc,$ and $id$ is the identity
matrix of $Sp_{2n}(\C).$ Note that $\C^{*}_I=\C_{J}^{*}$ if and
only if $\zeta^I$ and $\zeta^J$ are in the same orbit. Similarly,
 $\R^{*}_I=\R_{J}^{*}$ if and
only if $\zeta^I$ and $\zeta^J$ are in the same orbit. Therefore
there are $|O_n|$-complex curves and $|o_n|$-real curves in
$\mcvc,$ both ramified at the identity matrix $id \in Sp_{2n}(\C)$
$(t=0).$ The situation with $\mcvd$ is exactly same. The same
arguments also applies to the cases $\mcvb,$ and so there are
$|O_{n+1}|$-complex curves and $|o_{n+1}|$-real curves in $\mcvb,$
both ramified at the identity matrix $id \in SO_{2n+1}(\C)$
$(t=0).$

\section{EVALUATIONS OF SYMMETRIC FUNCTIONS AT ROOTS OF UNITY}\label{sec:6}

In this section, we set up the orthogonality formulas of $\Qt$-
and  $\Pt$-polynomials at roots of unity,  complete the theorems
\ref{thm;ogo}, \ref{2nd} and \ref{thm:oge}, using these
orthogonality formulas, and determine on which $n$-tuples of roots
of unity all $\Pt$-polynomials are positive.

\subsection{Orthogonality formulas for $\Qt$- and $\Pt$-polynomials}
For $I\in \mathcal{T}_n,$ there is $n$-tuple
$(\hat{i}_1,...,\hat{i}_n) \in \mathcal{T}_n$ such that the two
sets $\{\zeta^{i_1},...\zeta^{i_n}\}$ and
$\{\zeta^{\hat{i}_1},...,\zeta^{\hat{i}_n}\}$ enumerate all roots
of $(-1)^{n+1}.$ Denote $\hat{J}=(\hat{i}_1,...,\hat{i}_n).$

\begin{lemma}[\cite{Riet1}]\label{rie}
Let  $E_k$ and $H_k$ the elementary and complete symmetric
polynomials, respectively. Then for $I, J\in \mathcal{I}_n,$ we
have the following identities, \bn

\item
 $\prod_{k=1}^{n}
\prod_{l=1}^{n}(1-\zeta^{i_k}\zeta^{-\hat{j}_l})=\delta_{I,J}
\displaystyle{\frac{(2n)^n}{|\textrm{Vand}(\zeta^I)|^2}},$
\item $E_k(\zeta^{-I})=H_k(\zeta^{I^{*}}),$ \en
where
$\textrm{Vand}(\zeta^J):=\prod_{k<l}(\zeta^{j_k}-\zeta^{j_l}).$
\end{lemma}

\bp For proof of these identities, we refer to Proposition $4.3$
and Lemma $4.4$ in \cite{Riet1}.\ep

The next Proposition $\ref{pro7.2.1}$ are the analogue for $\Qt$-
and $\Pt$-polynomials of Propositions $4.3$ and $6.1$ in
\cite{Riet1}.
\begin{proposition}\label{p}
Let $I, J \in \mathcal{I}_n$ and $z_1,...,z_n, t\in \C$. Then we
have \bn\label{pro7.2.1}
 \item $\sum_{\lambda\in
\mathcal{R}(n)} {P_\lambda(\zeta^{J^*})}
{\Qt_\lambda(z_1,...,z_n)}=\prod_{k=1}^{n}\prod_{l=1}^{n}(1-z_k{\zeta^{-\hat{j}_l}}),$

 \item $\sum_{\lambda\in \mathcal{R}(n)}{P_\lambda(\zeta^{J^\ast})}{\Qt_\lambda(\zeta^I)}=\delta_{I,J}\displaystyle{\frac{(2n)^n}
 {|\textrm{Vand}(\zeta^I)|^2}},$
 \item $\sum_{\lambda\in
\mathcal{R}(n)} \overline{\omega(P_{\lambda})(\zeta^J)}
{\widetilde{Q}_\lambda(\zeta^I)}=\delta_{I,J}\displaystyle{\frac{(2n)^n}{|\textrm{Vand}(\zeta^I)|^2}},$
\item $\sum_{\lambda \in \mathcal{D}(n)}\Pt_{\lambda}(\zeta^I)\Pt_{
\hat{\lambda}}(\zeta^J)=\delta_{I,J}S_{\rho_n}(\zeta^I),$\en
 where
$\omega$ is the involution on the ring of symmetric polynomials
interchanging the elementary and the complete symmetric
polynomials, $\bar{()}$ stands for complex conjugation, and
$\rho_n:=(n,n-1,...,1)$.

\end{proposition}

\begin{proof}
We have the following identity from  p.234 in \cite{Mac1},

\begin{equation}\label{mac}
 \sum_{\lambda\in
\mathcal{R}(n)}{P_\lambda(z_1,...,z_n)}Q_\lambda^\prime(w_1,...,w_n)=
\prod_{k=1}^{n}\prod_{l=1}^{n}\frac{1}{(1-z_k{w_l})},
\end{equation}
where  $Q_{\lambda}^{\prime}$ are the Hall-Littlewood functions
 (Proposition $4.9$ in \cite{Pra and Rat 1}).
Note that from the definition of $Q_{\lambda}^{\prime},$ we have
$\omega(Q^{\prime}_\lambda)=\Qt_\lambda.$ Therefore if we take the
involution $\omega$ with respect to $W=(w_1,...,w_n)$ on the both
sides of (\ref{mac}), we have

\begin{displaymath}
\sum_{\lambda\in
\rn}{P_\lambda(z_1,...,z_n)}{\widetilde{Q}_\lambda(w_1,...,w_n)}=\prod_{k=1}^{n}\prod_{l=1}^{n}(1+{z_k}w_l).
\end{displaymath}
Now we replace $(w_1,...,w_n)$ by $\zeta^{J^*}$, and since
$\zeta^n=-1,$ we get $(1).$\\
For $(2),$  replace $(z_1,...,z_n)$ by
$(\zeta^{i_1},...,\zeta^{i_n})$. Then the formula comes from Lemma
$\ref{rie}$ $(1)$. The formula $(3)$  is a direct consequence of
Lemma $\ref{rie}$ $(2).$ For $(4),$ consider the following
identity from Lemma $2.7$ in \cite{Las and Pra},

\begin{displaymath}\sum_{\lambda\in \mathcal{D}(n)}\Pt_{\lambda}(X_n^w)\Pt_{\hat{ \lambda}}(X_n)=\left\{\begin{array}{cc}
S_{\rho_n}(X_n)& \textrm{if} \hspace{0.05in} w \in S_n,\\
0  & \textrm{if} \hspace{0.05in} w \in W_n\setminus S_n,
\end{array}\right.
\end{displaymath}
where  $w$ acts on $(x_1,...,x_n)$ naturally, i.e., $s_i$
interchanges $x_i$ and $x_{i+1},$ and $\dot{s}_n$ interchanges
$\pm x_n$ and $\mp x_{n+1}.$ The formula $(4)$ follows directly
from the fact that for given $I,J \in \mathcal{I}_n$ with $I\ne
J,$ there is $w\in W_n \setminus S_n$ such that
$(\zeta^{i_1},...,\zeta^{i_n})^w=(\zeta^{j_1},...,\zeta^{j_n}).$
\end{proof}

\begin{remark}
Let  $A=[A_{I,\lambda}]$  be a square matrix defined by
$A_{I,\lambda}=\Pt_{\lambda}(\zeta^{I})$ for $I\in \mathcal{I}_n$
and $\lambda \in \mathcal{D}(n).$ Then it follows from $(4)$ of
Proposition \ref{p} that the row vectors of $A$ are linearly
independent, and so there is no zero row  vector or zero column
vector in $A.$ That is, there are nonzero entries in all rows and
columns. Furthermore we have another kind of orthogonality
properties of $\Qt$-polynomials (and hence $\Pt$-polynomials)
given as follows.
\end{remark}

\begin{corollary}\label{corollary column ortho}
For any $\lambda, \mu\in \mathcal{D}(n),$ we have the following
identities

\bn
 \item $\displaystyle{\sum_{I\in
\mathcal{I}_n}\Pt_{\lambda}(\zeta^I)\Pt_{
\hat{\mu}}(\zeta^I)\displaystyle{\frac{1}{S_{\rho_n}(\zeta^I)}}}=\delta_{\lambda,\mu},$
\item
$\displaystyle{\frac{1}{(2n)^n}}\sum_{I\in
\mathcal{I}_n}{P_\lambda(\zeta^{I^\ast})}{\Qt_\mu(\zeta^I)}\displaystyle{
 {|\textrm{Vand}(\zeta^I)|^2}}=\delta_{\lambda,\mu},$
 \item $\displaystyle{\frac{1}{(2n)^n}}\sum_{I\in
\mathcal{I}_n}{\overline{\omega(P_\lambda)(\zeta^{I^\ast})}}{\Qt_\mu(\zeta^I)}\displaystyle{
 {|\textrm{Vand}(\zeta^I)|^2}}=\delta_{\lambda,\mu}.$
\en
\end{corollary}

\begin{proof}
For $(1),$ consider the matrix $B=[B_{\mu,J}]$ defined by
$$B_{\mu,J}:=\frac{\Pt_{\hat{\mu}}(\zeta^J)}{S_{\rho_n}(\zeta^J)}$$
for $\mu\in \mathcal{D}(n)$ and $J\in \mathcal{I}_n.$ Then the
entries of the multiplication $AB$ are given by the left hand side
of $(4)$ of Proposition \ref{p}. Hence $AB$ is equal to the
identity matrix, and so is the multiplication
$BA:=[C_{\mu,\lambda}].$ The identity $(1)$ follows since the left
hand side of $(1)$ is $C_{\mu,\lambda}.$ The identities $(2)$ and
$(3)$ follow by applying the same idea.
\end{proof}

\subsection{Completion of proofs of theorems \ref{thm;ogo}, \ref{2nd} and \ref{thm:oge}.}\label{completion of proofs}
 To complete the proof of Theorem \ref{thm;ogo}, we have to
show that $\phi_C$ is injective. Let $F$ be element of
$\mathfrak{R}_o$ such that $\dot{F}=0,$ i.e., $\dot{F}(u)=0$ for
all $u \in \mcvc.$ Applying the quantum Pieri rule for $\ogo$(
Corollary $5$ of \cite{KT1}) repeatedly, $F$ can be written as
$F=\sum_{\lambda \in \dn}F_{\lambda}\Pt_\lambda,$ where $F_\lambda
\in \C[q].$ Then for any $t\in \C$ and $I\in \In,$ we have
$$0=\dot{F}(u_{2n}(t\zeta^I))=\sum_{\lambda \in
\dn}F_{\lambda}(u_{2n}(t\zeta^I))\dot{\Pt}_\lambda(u_{2n}(t\zeta^I))=0,$$
 and hence
$$0=F(t\zeta^I)=\sum_{\lambda \in
\dn}F_{\lambda}(t^{2n})\Pt_\lambda(t\zeta^I)=\sum_{\lambda \in
\dn}t^{l(\lambda)}F_{\lambda}(t^{2n})\Pt_\lambda(\zeta^I),$$ where
$F_\lambda(t^{2n})$ is  a polynomial in $t^{2n}.$  Therefore for
fixed $t,$ the vector
$(t^{l(\lambda)}F_{\lambda}(t^{2n}))_\lambda$ lies in the
orthogonal component of the subspace of $\C^{|\dn|}$ spanned by
the vectors  $(\Pt_{\lambda}(\zeta^I))_\lambda$ indexed by $I\in
\In.$ But from Remark subsequent to Proposition \ref{p}, the
vectors $(\Pt_\lambda(\zeta^I))_{\lambda \in \dn}$ indexed by
$I\in \In$ are linearly independent, and since $|\dn|=|\In|,$ they
form a basis for $\C^{|\dn|}.$ Therefore the vector
$(t^{l(\lambda)}F_{\lambda}(t^{2n}))_\lambda$ is a zero vector.
Since this holds for all $t\in \C,$ $F_{\lambda}$ is $0$ in
$\C[q]$ for all $\lambda \in \dn.$ Therefore $\phi_C$ is
injective. The same argument applies to the cases and $\oge$ and
$\lg$ to complete the proof of the theorems $\ref{thm:oge}$ and
$\ref{2nd}.$

\begin{corollary}\label{basis 1}
Let $X_{i,i}$ be given as in $(\ref{2 paffian}),$ $\mathfrak{I}_C$
the ideal of $\Z[X_1,...,X_n]$ generated by $X_{i,i},$
$i=1,...n-1,$  and let
$\mathfrak{R}^{\Z}_C:=\Z[X_1,...,X_n]/\mathfrak{I}_C.$ Then the
set $\{\dot{q}\dot{\Pt}_\lambda \mid k\in \Z_{\geq 0}, \lambda \in
\dn\}$ forms a basis for the ring $\mathfrak{R}^\Z_C.$
\end{corollary}
\begin{proof}
By Theorem \ref{thm;ogo} and the Quantum Pieri Rule for $\ogo$ (
Corollary $5$ of \cite{KT1}), any element of $\mathfrak{R}^\Z_C$
can be written as a $\Z$-linear combination of
$\dot{q}^k\dot{\Pt}_\lambda.$ By the same argument as above, they
are linearly independent.
\end{proof}

\begin{corollary}\label{basis 2}
Let $X_{i,i}$ be given as in $(\ref{2 paffian}),$ $\mathfrak{I}_B$
the ideal of $\Z[X_1,...,X_{n+1}]$ generated by $X_{i,i},$
$i=1,...,n,$ and let
$\mathfrak{R}^\Z_B:=\Z[X_1,...,X_{n+1}]/\mathfrak{I}_B.$ The set
$\{\dot{q}\dot{\Qt}_\lambda \mid k\in \Z_{\geq 0}, \lambda \in
\dn\}$ forms a basis for the ring $\mathfrak{R}^\Z_B.$
\end{corollary}
\begin{proof}
The same argument as in the proof of Corollary \ref{basis 1}
applies.
\end{proof}

\subsection{Positivity of $\Qt$- and $\Pt$-polynomials}
\begin{lemma}[Perron-Frobenius Theorem]
Suppose a matrix $A$ is nonnegative, i.e., all entries $A_{i,j}$
are nonnegative. Then there is a nonnegative eigenvalue $K$ with a
nonnegative eigenvector and such that $K$ is maximal among
absolute values of all eigenvalues of $A.$
\end{lemma}
\begin{proof}
We refer to \cite{Minc} for a proof.
\end{proof}

\begin{lemma}\label{lemma;nonzero}
For any $I\in \mathcal{I}_n,$ we have $\Pt_{\rho_n}(\zeta^I)\ne0.$
\end{lemma}
\begin{proof}
From Corollary $7$ of \cite{KT1}, we have that
$\tau_{\rho_n}\cdotp \tau_{\rho_n}=\tau_n q^{n/2}$ if $n$ is even,
and $\tau_{\rho_n}\cdotp \tau_{\rho_n}=q^{(n+1)/2}$ if $n$ is odd.
From the identification of $\tau_\lambda$ and $\Pt_\lambda$, and
the fact (Lemma \ref{lemma:quant12}) that $\Pt_n(\zeta^I)$ is
nonzero for any $I\in \mathcal{I}_n,$  it follows that
$\Pt_{\rho_n}(\zeta^I)$ is nonzero.
\end{proof}
\begin{theorem}\label{positivity of Q}
We have the following inequalities.
 \bn \item
$\Pt_{\lambda}(\zeta^{I_0})>0$ for all strict partitions $\lambda
\in \dn.$ And  for $I^s \in \In$ with $I\ne I_0,$ there exists
strict partitions $\mu \in \dn$
 such that $\Pt_\mu(\zeta^I)<0.$
 Furthermore, if $I\in \In$ is such that  $\Pt_\lambda(\zeta^I)\geq 0$ for all
 $\lambda \in \dn,$ then $I=I_0.$
\item For any $\lambda\in \dn,$ $\Pt_{\lambda}(\zeta^{I_0})\geq
|\Pt_{\lambda}(\zeta^{I})|$ for all $I \in \In.$ \en
\end{theorem}

\begin{proof}
The idea of the proof will be taken more or less verbatim from
Proposition $11.1$ of \cite{Riet1}, where the analogue of $(2)$
was obtained. Consider
$qH^*_\C(OG^o(n))_{q=1}:=qH^*_{\C}(OG^o(n))\otimes \C[q]/(q-1),$
the specialization of the quantum cohomology ring at $q=1.$ By the
isomorphism of Theorem \ref{thm;ogo},  this ring may be viewed as
the coordinate ring of the zero-dimensional subvariety
$\{u_{2n}(\zeta^I)\mid I \in \In\}$ of $\mcvc$ cut out by the
function $\dot{q}-1 \in \mo(\mcvc)$ and has a basis given by the
Schubert basis $\tau_{\lambda},$ $\lambda \in \dn$. For  $I \in
\In,$ define $\tau_I:=\sum_{\nu \in
\dn}\Pt_\nu(\zeta^I)\tau_{\hat{\nu}}.$ Then, by Proposition
$\ref{p}$ $(4),$ we have that
$$\tau_I(u_{2n}(\zeta^J))=\delta_{I,J}S_{\rho_n}(\zeta^J).$$
Therefore  the vectors $\tau_I,$ $I\in \In,$ are linearly
independent, and hence form a basis for $qH^*_\C(OG^o(n))_{q=1}.$
Define the multiplication operator $[\tau_\lambda]:\tau \mapsto
\tau_{\lambda}\cdot \tau$ on $qH^*_\C(OG^o(n))_{q=1}$.  Then for
each $I \in \In,$ we have
$[\tau_\lambda](\tau_I)=\Pt_\lambda(\zeta^I)\tau_I$ for all $I \in
\In.$ So  $[\tau_\lambda]$ has  eigenvalues
$\Pt_\lambda(\zeta^I),$ $I\in \In,$ with eigenvectors $\tau_I.$ In
fact the set $\{\tau_I |I \in \tI\}$ is a simultaneous eigenbasis
for the operators $[\tau]$ with  $\tau \in
qH^*_\C(OG^o(n))_{q=1}$. Let $S$ be the set of elements $\tau$ in
$qH^*_\R(OG^o(n))_{q=1}$ such that $[\tau]$ has simple
eigenvalues. Then $S$ is open dense subset of
$qH^*_\R(OG^o(n))_{q=1},$ and so all $[\tau_\lambda]$ can be
approximated to arbitrary precision by the operators $[\tau]$ with
$\tau \in S$; for any small $\epsilon
>0, $ there are $0<\epsilon_{\mu}<\epsilon$ such that $[\tau_\lambda]+\sum_{\mu \in
\dn}\epsilon_{\mu}[\tau_\mu]$ has simple eigenvalues.

 Let $A_{\lambda}$ be the matrix of the operator $[\tau_\lambda]$ with
respect to the basis $\{ \tau_\nu \mid \nu \in \dn\}.$  Since the
entries $A_{\mu,\nu}^\lambda$ are Gromov-Witten invariants,
$A_\lambda$ is nonnegative and so is the matrix
$A_{\lambda}^\epsilon:=A_{\lambda}+\sum_{\mu}\epsilon_{\mu}
A_{\mu}.$ We also note that the vector $\tau_I$ has the
coordinates $(\Pt_{\lambda}(\zeta^I))_{\lambda\in \dn }$ with
respect to the basis $\{ \tau_\nu \mid \nu \in \dn\}.$ First
consider the case $\lambda=(1).$  In this case $A_{(1)}$ has
$|\In|$-distinct eigenvalues $\Pt_1 (\zeta^I),$ $I\in \In,$ with
eigenvectors $(\Pt_\lambda(\zeta^I))_{\lambda \in \dn}$. Therefore
all the eigenvalues of $A_{(1)}$ are simple, and, by the
Perron-Frobenius theorem,  we have the maximum eigenvalue
$\Pt_{1}(\zeta^{I_0})=\frac{1}{2}E_1(\zeta^{I_0})$ with a
nonnegative eigenvector, which is
$(\Pt_{\lambda}(\zeta^{I_0})_{\lambda \in \dn}$
  or
$-(\Pt_{\lambda}(\zeta^{I_0})_{\lambda \in \dn}$. But the vector
$-(\Pt_{\lambda}(\zeta^{I_0})_{\lambda \in \dn}$ can not be
nonnegative since $\Pt_{1}(\zeta^{I_0})> 0.$ Therefore
$(\Pt_{\lambda}(\zeta^{I_0})_{\lambda \in \dn}$ is a nonnegative
eigenvector, and so  $\Pt_{\lambda}(\zeta^{I_0})\geq 0$ for all
$\lambda \in \dn$.  For the second statement of $(1)$, fix $I\in
\mathcal{I}_n$ with $I\ne I_0.$ Since $\Pt_{\rho_n}(\zeta^I)$ is
no zero by Lemma \ref{lemma;nonzero}, Proposition \ref{p} $(4)$
implies that there is a strict partition $\mu\in\mathcal{D}(n)$
such that $\Pt_{\mu}(\zeta^I)<0.$  This proves the second
statement of $(1).$ The last statement of $(1)$ follows from the
second statement of $(1)$ and the fact that $I\in \In^s$ if and
only if $\Pt_{\lambda}(\zeta^I)$ are real for all $\lambda \in
\dn.$ This proves $(1),$ except that $\Pt_\lambda(\zeta^{I_0})$
are strictly positive for all $\lambda\in \mathcal{D}(n)$. This
will be done after proving $(2).$ For $(2),$ we apply the
Perron-Frobinius Theorem above to $A_{\lambda}^\epsilon$. Then
 there is $I\in \In$ such that  $$K_\lambda^\epsilon(I):=\Pt_\lambda(\zeta^I)+
 \sum_{\mu} \epsilon_\mu \Pt_{\mu}(\zeta^I)$$ is a maximal real eigenvalue with a
 nonnegative eigenvector, and such that
 $K_\lambda^\epsilon(I)\geq |K_\lambda^\epsilon(J)|$ for all $J\in \In,$
But since the vector $(\Pt_{\mu}(\zeta^{I_0}))_{\mu \in \dn}$
 is a unique nonnegative eigenvector(up to positive constants) by the last statement of $(1)$, we
have $I=I_0,$  and for all $I \in \In,$ $$K_\lambda^\epsilon(I_0)
\geq |K_\lambda^\epsilon(I)|.$$ Since this is true for all
$\epsilon >0,$ we have
$\Pt_\lambda(\zeta^{I_0})\geq|\Pt_\lambda(\zeta^{I})|$ for all $I
\in \In.$ This proves $(2).$ Now we prove that
$\Pt_\lambda(\zeta^{I_0})$  are strictly positive for all
$\lambda\in \mathcal{D}(n)$. If $\Pt_\lambda(\zeta^{I_0})$ is zero
for some $\lambda \in \mathcal{D}(n),$  by $(2)$,
$\Pt_\lambda(\zeta^{I})$ is zero for all $I\in \mathcal{I}_n$,
which is impossible by Remark subsequent to Proposition \ref{p}.
Therefore $\Pt_\lambda(\zeta^{I_0})$ is strictly positive for all
$\lambda \in \mathcal{D}(n).$ This completes the proof.
\end{proof}

\section{GROMOV-WITTEN INVARIANTS}\label{sec:7}

In this section, we use the orthogonality formula given in
 Section \ref{sec:6} and the theorems, $\ref{thm;ogo}$ and $\ref{2nd},$ to
derive the Vafa-Intriligator type formulas, which compute the
3-point, genus zero, Gromov-Witten invariants for $\ogo,$ ($\oge$)
and $\lg,$ and also as an application we give an analogue of
Poincar\'{e} duality pairing on $qH^*(\oge)$ and $qH^*(\oge)$. We
start with the following preliminary results which will be crucial
to find the Vafa-Intriligator type formulas.

\begin{lemma} \label{first}
Let  $I\in \In,$ $t\in\C^{\ast},$ $P$ be a homogeneous symmetric
polynomial in the variables $x_1,...,x_n$.Then we have the
following expressions,
$$
P(t\zeta^I)=\sum_{\nu\in
\mathcal{R}(n)}{m_\nu^P}(t)\Qt_\nu(t\zeta^I),
$$ where $m_\nu^P(t)=\frac{1}{(2n)^n}\sum_{J\in
\In}P(t\zeta^J)P_{\nu}(t^{-1}\zeta^{J^{*}})|
\textrm{Vand}(\zeta^J)|^2.$
\end{lemma}

\bp By Proposition \ref{p} $(2)$, we have the expression for
$P(t\zI)$

\begin{displaymath}
P(t\zI)=\frac{1}{(2n)^n}\sum_{J\in \tI }P(t\zJ)|
\textrm{Vand}(\zJ)|^2 \sum_{\nu \in \mathcal{R}(n)}
P_\nu(\zeta^{J^{\ast}})\Qt_\nu(\zI),
\end{displaymath}
which is equal to
\begin{displaymath}
\frac{1}{(2n)^n}\sum_{\nu \in \mathcal{R}(n)} \sum_{J\in \tI}
P(t\zJ)|\textrm{Vand}(\zJ)|^2P_\nu(t^{-1}\zeta^{J^{\ast}})\Qt_{\nu}(t\zeta^{I}).
\end{displaymath}
This proves the lemma. \ep

\begin{lemma}[\cite{Riet1}] \label{m}
Let $m$ be a homogeneous polynomial in n variables. Then we have\\
$\sum_{J\in \tI} m(\zJ)=0$  unless $\textrm{deg}(m)\equiv0
\hspace{0.08in} mod (2n).$ In particular,
$$m_\nu^P(t)=\frac{1}{(2n)^n}t^{(\textrm{deg} P- |\nu
|)}\sum_{J\in
\In}P(\zeta^J)P_{\nu}(\zeta^{J^{*}})|\textrm{Vand}(\zeta^J)|^2$$
if $\textrm{deg} P \equiv | \nu|$ mod $2n$, and otherwise
$m_\nu^P(t)=0$.
\end{lemma}

\bp Let  $ M(t):=\sum_{J\in \In} m(t \zJ).$ Then we have the
relation $M(t)=M(t \zeta),$ and hence $M(t)$ is a polynomial in
$t^{2n}.$ Therefore unless $\textrm{deg}$ $m$ is divisible by
$2n,$ we must have $M(t)=t^{\textrm{deg} m}\sum_{J\in \In} m(
\zJ)=0$  for all $t.$ It follows  that $ \sum_{J\in \In} m(t
\zJ)=0$ unless $\textrm{deg}$ $m$ is divisible by $2n.$ Apply the
same analysis to $m_\nu^P(t)$  to get the rest of the lemma.
 \ep
\subsection{Vafa-Intriligator type formula for $\ogo$ $(\oge)$}

\begin{proposition}\label{pro 8.1.1}
 Let $\dot{P}$ be an element of $\mathfrak{R}_C^\Z \subset \mo(\mcvc).$ Then $\dot{P}$ can be expressed as

\be \label{Pdot} \dot{P}=\sum_{\nu \in \dn, k\in
\Z_{\geq0}}\dot{P}_{\nu,k}\hspace {0.07in}{\dot{q}}^k
\hspace{0.07in}\dot{\Pt}_{\nu}.\ee

Here $\dot{P}_{\nu,k}$  is an integer defined as follows. If $\nu$
is a strict partition such that \\$\textrm{deg}
\hspace{0.07in}\dot{P}=|\nu| + 2nk $ for some nonnegative integer
$k,$

\be \label{general gro}
\dot{P}_{\nu,k}:=\frac{2^{l(\nu)+2k}}{(2n)^n}\sum_{m=0}^{a(\nu)}\sum_{J\in
\In}
P(\zeta^J){P_{\nu(m)}(\zeta^{J^{\ast}})}|\textrm{Vand}(\zJ)|^2
,\ee

\bigskip
 and otherwise  $\dot{P}_{\nu,k}:=0,$

where $a(\nu)=\lfloor \frac{n-l(\nu)}{2}\rfloor ,$ and
$\nu(m)=((n)^{(2m)},\nu_1,...,\nu_l)$.
\end{proposition}
\bp From Corollary \ref{basis 1}, the set $\{
\dot{\Pt}_{\nu}\dot{q}^k \mid k\in \Z_{\geq0}, \nu \in \dn\}$ is a
$\Z$-basis for  the ring $\mathfrak{R}_C^\Z$, so $\dot{P}_{\nu,k}$
are integers. To get $(\ref{general gro}),$ we consider the
following expression for $P(t\zeta^I)$ from Lemma $\ref{first},$

\be \label{initial1} P(t\zeta^I)=\sum_{\nu\in
\mathcal{R}(n)}\sum_{J\in
\In}\frac{|\textrm{Vand}(\zJ)|^2}{(2n)^n}P(\zJ)P_{\nu}(\zeta^{J^{\ast}})t^{(\textrm{deg}
P-|\nu|)} \Qt_{\nu}(t\zI).\ee

By Lemma $\ref{m},$ the terms in $(\ref{initial1})$ corresponding
to $\nu$ with $ 2n \nmid ( \textrm{deg} P-\mid \nu \mid) $ vanish,
and so $P(t\zeta^I)$ can be written as \be
P(t\zeta^I)=\sum_{\nu\in A}\sum_{J\in
\In}\frac{|\textrm{Vand}(\zJ)|^2}{(2n)^n}P(\zJ)P_{\nu}(\zeta^{J^{\ast}})t^{(\textrm{deg}
P-|\nu|)} \Qt_{\nu}(t\zI),\ee where the index set $A$ is the set
of all partitions $\nu \in \rn$ such that
$\textrm{deg}\hspace{0.03in} P-|\nu|$ is divisible by $2n.$

If we evaluate the identity of functions on $\mcvc$ determined by
this identity, then the terms  $\Qtd_{\nu}$ vanish for nonstrict
partitions $\nu$ with repeated parts from $\{1,2,...,n-1\},$ by
the fact that  $\Qtd_{i,i}=0 $ in the ring
$\mathcal{O}(\mathcal{V}_n^C)$ for $i=1,...,n-1,$   and by the
factorization property  of $\Qt$-polynomials ($ (\ref{fac:Q})$ in
$\ref{subsec:sym poly}$). Therefore, for
$u(t\zeta^I)=u_{2n}(t\zeta^I)\in \mcvc,$ we have
 $$
\dot{P}(u(t\zeta^I))=\sum_{\nu\in B}\sum_{J\in \In}\frac{|
\textrm{Vand}(\zJ)|^2}{(2n)^n}\dot{P}(u(\zJ))\dot{P}_{\nu}(u(\zeta^{J^{\ast}}))
t^{(\textrm{deg} P-|\nu|)}\Qtd_{\nu}(u(t\zI)),
$$
where $B$ is the subset of $A$ consisting $\nu$ such that $\nu$
has no repetitions from $\{1,...,n-1\}.$ Note that each partition
in $B$ is of the form $\nu(m)$ for some strict partition $\nu\in
\mathcal{D}(n)$ and a nonnegative integer $m.$  Using the
factorization property  of $\Qt$-polynomials ((\ref{fac:Q}) in
\ref{subsec:sym poly}), and the identities,
$\Qt_{(n^{2m})}={\Qt_{(n^2)}}^m,$ and
$\dot{\Qt}_{(n^2)}(u(t\zeta^I))=t^{2n},$  we can rewrite
$$ \dot{P}(u(t\zeta^I))=\sum_{\nu\in
\mathcal{D}(n)}\sum_{m=0}^{a(\nu)}\sum_{J\in \In}\frac{|
\textrm{Vand}(\zJ)|^2}{(2n)^n}\dot{P}(u(\zJ))\dot{P}_{\nu(m)}(u(\zeta^{J^{\ast}}))
t^{(\textrm{deg} P-|\nu(m)|)}\Qtd_{\nu(m)}(u(t\zI))$$

$$ =\sum_{\nu \in \mathcal{D}(n)}\sum_{m=0}^{a(\nu)}\sum_{J\in
\In}\frac{|\textrm{Vand}(\zJ)|^2}{(2n)^n}\dot{P}(u(\zJ))\dot{P}_{\nu(m)}(u(\zeta^{J^{*}}))t^{(\textrm{deg}P
-|\nu(m)|)} \Qtd_{\nu}(u(t\zI))\Qtd_{(n^{2m})}(u(t\zI)) $$

\begin{displaymath}
=\sum_{\nu \in \mathcal{D}(n)}\sum_{m=0}^{a(\nu)}\sum_{J\in
\In}\frac{|\textrm{Vand}(\zJ)|^2}{(2n)^n}P(\zeta^J)P_{\nu(m)}(\zeta^{J^{*}})
\Qtd_{\nu}(u(t\zI)){\dot{\Qt}_{(n^2)}(u(t\zI))}^{\frac{deg P- |\nu
|}{2n}}.
\end{displaymath}

Here the upper bound $a(\nu)$ is taken since $\nu(m)$  belong to
$\rn.$

Since $u(t\zeta^I)$ are arbitrary points in $\mcvc$,
$\dot{\Qt}_{\nu}=2^{l(\nu)}\dot{\Pt}_{\nu}$ and
$\dot{\Qt}_{(n^2)}=4\dot{q},$ the function $\dot{P}$ can be
written as in $(\ref{Pdot})$, with $\dot{P}_{\nu,k }$ defined as
in $(\ref{general gro}).$\ep

\begin{corollary}\label{corollary:gro;oge}

For partitions $\lambda,$ $\mu,$ $\nu \in \mathcal{D}(n),$ and a
nonnegative integer $k,$ the Gromov-Witten invariant
$<\tau_{\lambda},\tau_{\mu},\tau_{\hat{\nu}}>_k$ for $\ogo (\cong
\oge) $
 is given by
\be~\label{gro:og}
<\tau_{\lambda},\tau_{\mu},\tau_{\hat{\nu}}>_k=\frac{2^{l(\nu)+2k}}{(2n)^n}\sum_{m=0}^{a(\nu)}\sum_{J\in
\In}\widetilde{P}_{\lambda}(\zeta^{J})\widetilde{P}_{\mu}(\zeta^{J})P_{\nu(m)}(\zeta^{J^{*}})|\textrm{Vand}(\zeta^J)|^2,
\ee
 whenever $ |\lambda| + |\mu|= |\nu| +2nk$, and otherwise by
$<\tau_{\lambda},\tau_{\mu},\tau_{\hat{\nu}}>_k=0.$
\end{corollary}
\bp Apply $\dot{P}=\dot{\Pt}_\lambda \dot{\Pt}_\mu$ to Proposition
\ref{pro 8.1.1}. Then the Gromov-Witten invariant
$<\tau_{\lambda},\tau_{\mu},\tau_{\hat{\nu}}>_k $ is the
coefficient of $\dot{\Pt}_\nu \dot{q}^k$ in the expansion of
$\dot{\Pt}_\lambda \dot{\Pt}_\mu$ in the basis  $\{
\Pt_{\nu}\dot{q}^k \mid  \nu \in \dn \hs \textrm{and} \hs k\geq 0
\},$ which is given as above. \ep

\subsection{Vafa-Intriligator type formula for $\lg$}

\begin{proposition}\label{fifth}%%prop 10.1
Let $\dot{P}$ be an element of $\mathfrak{R}_B^\Z \subset
\mo(\mcvb).$ Then $\dot{P}$ can be expressed as
\begin{displaymath}
\dot{P}=\sum_{\nu\in \dn, d\in \Z_{\geq0}}\dot{P}_{\nu,d}\hspace
{0.07in}{\dot{q}}^d \hspace{0.07in}\Qtd_{\nu}.
\end{displaymath}

Here $\dot{P}_{\nu,d}$  is an integer defined as follows. If $\nu$
is a strict partition in $\dn$ with  $\textrm{deg}\hspace
{0.05in}\dot{P}=|\nu| + (n+1)d $ for some nonnegative integer $d,$

\be \label{vafa}
\dot{P}_{\nu,d}:=\frac{1}{2^d(2n+2)^{(n+1)}}\sum_{m=0}^{b(\nu)}\sum_{J\in
\Im}
P(\zeta^J){P_{\nu[m]}(\zeta^{J^{\ast}})|\textrm{Vand}(\zJ)|^2}
.\ee

 and otherwise  $\dot{P}_{\nu,d}:=0$,
where $\nu[m]$ and $b(\nu)$ are defined as
 \bdm \nu[m]:= \left
\{\begin{array}{cc}
((n+1)^{2m},\nu_1,...\nu_l)  \hs\textrm{if}\hs d\hs \textrm{is even},\\
((n+1)^{2m+1},\nu_1,...,\nu_l) \hs\textrm{if}\hs d \hs \textrm{is
odd},
 \end{array}\right. \edm

 \bdm b(\nu):= \left
\{\begin{array}{cc}
\lfloor\frac{n+1-l(\nu)}{2}\rfloor \hs\textrm{if}\hs d \hs \textrm{is even},\\
\lfloor\frac{n-l(\nu)}{2}\rfloor \hs\textrm{if}\hs d \hs
\textrm{is odd}.
 \end{array}\right. \edm

\end{proposition}

\bp

From Corollary \ref{basis 2}, the set  $\{\dot{\Qt}_\nu \dot{q}^d
\mid \nu \in \dn,  d\geq 0 \}$ is a $\Z$-basis for the ring
$\mathfrak{R}_B^\Z,$ so $\dot{P}_{\nu,d}$ are integers. For
$(\ref{vafa}),$ we consider the following expression for $\dot{P}$
from  Lemma \ref{first},

\bdm \dot{P}(u(t\zI))=\sum_{\nu\in \rnp}\sum_{J\in \Im}\frac{|
\textrm{Vand}(\zJ)|^2}{(2n+2)^{(n+1)}}\dot{P}(u(t\zJ))\dot{P}_{\nu}(u(t^{-1}\zeta^{J^{\ast}}))
\Qtd_{\nu}(u(t\zI)), \edm where $u(t\zeta^I)=u_{2n+1}(t\zeta^I)\in
\mcvb. $

If we apply Lemma \ref{m} to the above identity,  the terms
$\Qtd_{\nu}$ vanish for the partitions $\nu$ with $ (2n+2) \nmid (
\textrm{deg} P- |\nu|), $ and if we
 evaluate  $\dot{P}$ on the points of $\mcvb$, then the terms
$\Qtd_{\nu}$  vanish for  nonstrict partitions $\nu$ with repeated
parts from $\{1,2,...,n\},$ since the variety  $\mcvb$ is defined
by the equations $\dot{\Qt}_{i,i}=0$ for $i=1,...,n.$ Therefore
$\dot{P}(t\zeta^I)$ is reduced to the following expression

\bdm\sum_{\nu \in A}\sum_{J\in \Im}\frac{|
\textrm{Vand}(\zJ)|^2}{(2n+2)^{n+1}}\dot{P}(u(t\zJ))\dot{P}_{\nu}(u(t^{-1}\zeta^{J^{\ast}}))
\Qtd_{\nu}(u(t\zI)), \edm where the index set $A$  is the set of
all partitions $\nu=(\nu_1,...,\nu_l) \in \rnp$  such that  if $
\nu_i\leq n$ for some $ i$, then $ \nu_{i}>\nu_{i+1}>\cdots >
\nu_l,$ and such that
 $\textrm{deg}P -|\nu|$ is divisible by $2(n+1).$
By applying the factorization property  of $\Qt$-polynomials
((\ref{fac:Qspe}) in \ref{subsec:sym poly}) and regrouping the
terms, we can write $\dot{P}(t\zeta^I)$ as follows
\begin{displaymath}
\sum_{\nu \in B}\sum_{m=0}^{b(\nu)}\sum_{J\in
\Im}\frac{|\textrm{Vand}(\zeta^J)|^2}{(2n+2)^{n+1}}\dot{P}(u(\zJ))\dot{P}_{\nu[m]}(u(\zeta^{J^{*}}))t^{(\textrm{deg}
P-|\nu(m)|)} \Qtd_{\nu}(u(t\zI))\Qtd_{((n+1)^{ r})}(u(t\zI)),
\end{displaymath}
where the index set $B$ is the set of all partitions $\nu \in \dn$
such that $\textrm{deg}P -|\nu|$ is divisible by $(n+1),$ $r=2m$
if $\frac{\textrm{deg}P-|\nu|}{n+1}$ is even, $r=2m+1$ otherwise,
and the second sum is taken up to $b(\nu)$ since $\nu[m]$ belongs
to $\rnp.$

Using the relation $q=2\Qt_{n+1},$ we can rewrite
$\dot{P}(t\zeta^I)$ as follows,
 \bdm \sum_{\nu\in
B}\sum_{m=0}^{b(\nu)}\sum_{J\in \Im}\frac{ |
\textrm{Vand}(\zeta^J)|^2}{k(P,n,\nu)}\dot{P}(u(\zJ))\dot{P}_{\nu[m]}(u(\zeta^{J^{*}}))
\Qtd_{\nu}(u(t\zI))\dot{q}^{\frac{\textrm{deg} P-
|\nu|}{n+1}}(u(t\zI)), \edm where
$k(P,n,\nu)=2^{\frac{\textrm{deg}P-|\nu|}{n+1}}(2n+2)^{n+1}.$

 If $d$ is a nonnegative integer such that $\textrm{deg}\hspace
{0.05in}\dot{P}= |\nu| + (n+1)d,$ and so $\frac{\textrm{deg} P-
|\nu|}{n+1}=d,$  the coefficient of
$\dot{\Qt_{\nu}}(u(t\zI))\dot{q}^{k}(u(t\zI))$ is
$\dot{P}_{\nu,d},$
 given as in $(\ref{vafa})$ .\ep

\begin{corollary} %% coro 10.2
For partitions $\lambda,\mu,$ and $\nu\in \dn,$ and a nonnegative
integer $d,$ the Gromov-Witten invariant
$<\sigma_{\lambda},\sigma_{\mu},\sigma_{\hat{\nu}}>_d$ for $\lg$
is given by \be\label{gromov:lg}
<\sigma_{\lambda},\sigma_{\mu},\sigma_{\hat{\nu}}>_d=\frac{1}{2^d(2n+2)^{(n+1)}}\sum_{m=0}^{b(\nu)}
 \sum_{J\in
 \Im}\Qtl(\zeta^{J})\Qt_{\mu}(\zeta^{J})P_{\nu[m]}(\zeta^{J^{*}})|
 \textrm{Vand}(\zeta^J)|^2,
\ee  whenever $|\lambda| + |\mu|=| \nu| +(n+1)d$, and otherwise by
$<\sigma_{\lambda},\sigma_{\mu},\sigma_{\hat{\nu}}>_d=0.$
\end{corollary}

\bp  If we apply $\dot{P}=\dot{\Qt}_\lambda \dot{\Qt}_\mu$ to
Proposition \ref{fifth}, the right hand side of
$(\ref{gromov:lg})$ is the coefficient of $\dot{\Pt}_\nu
\dot{q}^k$ in the expansion of $\dot{\Pt}_\lambda \dot{\Pt}_\mu$
in the basis $\{ \Pt_{\nu}\dot{q}^d \mid \nu \in \dn \hs
\textrm{and} \hs d \geq 0\},$ which is the Gromov-Witten invariant
$<\sigma_{\lambda},\sigma_{\mu},\sigma_{\hat{\nu}}>_d. $ \ep

For a partition $\nu \in \mathcal{D}(n-1) \subset \dn ,$ let
$\hat{\lambda}^n$ and $\hat{\lambda}^{(n-1)}$ be the partitions
that complement $\lambda$ in the sets $\{1,2,...,n\}$ and
$\{1,2,...,n-1\},$ respectively,  and denote
$\tilde{\nu}:=(n,\nu).$ We can relate the quantum cohomology rings
$\qhoge$ and $qH^*(LG(n-1))$ to each other in the following way

\begin{corollary}[\cite{KT1}]
Let $\lambda,$ $\mu$ and $\nu$ be partitions in
$\mathcal{D}(n-1).$ Then we have the following identities.
\be<\tau_{\lambda},\tau_{\mu},
\tau_{\hat{\nu}^n}>_k=2^{4k+l(\nu)-l(\lambda)-l(\mu)}<\sigma_\lambda,
\sigma_\mu, \sigma_{\hat{\nu}^{(n-1)}}>_{2k},\ee
\be\label{compare}
<\tau_\lambda,\tau_\mu,\tau_{\hat{\nu}^{(n-1)}}>_k=2^{4k+1+l(\tilde{\nu})-l(\lambda)-l(\mu)}
<\sigma_\lambda,\sigma_\mu,\sigma_{\hat{\nu}^{(n-1)}}>_{2k+1}.\ee
\end{corollary}
\bp These identities follow easily by comparing (\ref{gro:og}) and
(\ref{gromov:lg}). \ep
\subsection{Poincar\'{e} Duality for $qH^{*}(\oge)$ and
$qH^{*}(\lg)$}  Consider the Poincar\'{e} duality pairing
$$(\hspace{0.07in},\hspace{0.05in}):H^*(\oge)\times H^*(\oge)\rightarrow \Z$$ defined by $(\tau, \tau^{\prime})=(\tau \cup \tau^{\prime})[\oge],$
where $[\oge]$ is the fundamental class of $\oge.$ Recall that in
the Schubert basis $\{\tau_\lambda\mid\lambda \in
\mathcal{D}(n)\}$ for $\oge$ the pairing
$(\hspace{0.07in},\hspace{0.05in})$ is given by
$$(\tau_\mu,\tau_{\hat{\nu}})=\delta_{\mu,\nu}.$$
Now we give an analogue of this pairing for $qH^{*}(\oge).$ Let
$(\hspace{0.07in},\hspace{0.05in})_q$ be a $\Z[q]$-linear map
$$(\hspace{0.07in},\hspace{0.05in})_q:qH^{*}(\oge)\times qH^{*}(\oge)\rightarrow\Z[q]$$
which send $(\tau,\tau^{\prime})$ to the coefficient of
$\tau_{\rho_n}$ in the multiplication $\tau\cdotp\tau^{\prime}.$
Then the pairing $(\hspace{0.07in},\hspace{0.05in})_q$ specializes
to the classical Poincar\'{e} $(\hspace{0.07in},\hspace{0.05in})$
by setting $q=0.$ At the Schubert basis elements we have the
following evaluations; for $\mu, \nu \in \mathcal{D}(n),$ we have
$$(\tau_{\mu},\tau_{\hat{\nu}})_q=<\tau_{\mu},\tau_{\hat{\nu}},\tau_{\hat{\rho}_n}>_k q^k$$
if $|\mu|=|\nu|+2nk,$  otherwise
$(\tau_{\mu},\tau_{\hat{\nu}})_q=0$. By Corollary
\ref{corollary:gro;oge}, \be\label{formula100}
(\tau_{\mu},\tau_{\hat{\nu}})_q=\frac{2^{l(\nu)+2k}}{(2n)^n}q^k\sum_{m=0}^{a(\nu)}\sum_{J\in
\In}\widetilde{P}_{\mu}(\zeta^{J})P_{\nu(m)}(\zeta^{J^{*}})|\textrm{Vand}(\zeta^J)|^2\ee
if $|\mu|=|\nu|+2nk,$  otherwise
$(\tau_{\mu},\tau_{\hat{\nu}})_q=0$.
 From  the proof of Proposition \ref{pro 8.1.1} we see that the right hand side
 of (\ref{formula100}) is equal to
$$\frac{2^{l(\nu)}}{(2n)^n}t^{(|\mu|-|\nu|)}\sum_{m=0}^{a(\nu)}\sum_{J\in
\In}\widetilde{P}_{\mu}(\zeta^{J})P_{\nu(m)}(\zeta^{J^{*}})|\textrm{Vand}(\zeta^J)|^2.$$

Therefore we get the following result.
\begin{corollary} For any $\mu,\nu\in \mathcal{D}(n),$ the pairing $(\hspace{0.07in},\hspace{0.05in})_q$  has the same form
as the classical Poincar\'{e} duality pairing
$(\hspace{0.07in},\hspace{0.05in}),$ i.e.,
$$(\tau_{\mu},\tau_{\hat{\nu}})_q=\delta_{\mu,\nu}.$$
\end{corollary}
\begin{proof}
Note that since $\mu$ and $\nu$ are strict partitions, $\mu$ and
$\nu(m)$ coincide with each other precisely when $\mu=\nu$ and
$m=0.$ Therefore the corollary  follows directly from $(2)$ of
Corollary \ref{corollary column ortho}.
\end{proof}

 Similarly, we have the Poincar\'{e} duality pairing for
$H^{*}(\lg)$
$$(\hspace{0.07in},\hspace{0.05in}):H^*(\lg)\times H^*(\lg)\rightarrow \Z$$ defined by $(\sigma, \sigma^{\prime})=(\sigma \cup \sigma^{\prime})[\lg],$
where $[\lg]$ is the fundamental class of $\lg.$ Note that the
pairing $(\hspace{0.07in},\hspace{0.05in})$ satisfies the
Poincar\'{e} duality, that is, for any  $ \mu, \nu \in
\mathcal{D}(n),$ $(\sigma_\mu,
\sigma_{\hat{\nu}})=\delta_{\mu,\nu}.$ Now let
$(\hspace{0.07in},\hspace{0.05in})_q$ be a $\Z[q]$-linear map
$$(\hspace{0.07in},\hspace{0.05in})_q:qH^{*}(\lg)\times qH^{*}(\lg)\rightarrow\Z[q]$$
which send $(\sigma,\sigma^{\prime})$ to the coefficient of
$\sigma_{\rho_n}$ in the multiplication
$\sigma\cdotp\sigma^{\prime}.$ Then for any Schubert basis
elements $\sigma_\mu, \sigma_{\hat{\nu}}$, \be\label{fformula}
(\sigma_{\mu},\sigma_{\hat{\nu}})_q=<\sigma_{\mu},\sigma_{\hat{\nu}},\sigma_{\hat{\rho}_n}>_d
q^d\ee
 if $|\mu|=|\nu|+(n+1)d,$  otherwise
$(\sigma_{\mu},\sigma_{\hat{\nu}})_q=0$. By the same reasoning as
in the case of $\oge,$ the right hand side of $(\ref{fformula})$
is equal to
$$\frac{1}{(2n+2)^{(n+1)}}t^{(|\mu|-|\nu[m]|)}\Qt_{((n+1)^r)}(t\zeta^I)\sum_{m=0}^{b(\nu)}\sum_{J\in
\mathcal{I}_{n+1}}\widetilde{Q}_{\mu}(\zeta^{J})P_{\nu[m]}(\zeta^{J^{*}})|\textrm{Vand}(\zeta^J)|^2,$$
where $r=2m$ if $d$ is even, and $r=2m+1$ if $d$ is odd. Therefore
by applying $(2)$ of Corollary \ref{corollary column ortho}, we
get the following result.
 \begin{corollary} For any $\mu,\nu\in \mathcal{D}(n),$ the pairing $(\hspace{0.07in},\hspace{0.05in})_q$  has the same form
as the classical Poincar\'{e} duality pairing
$(\hspace{0.07in},\hspace{0.05in}),$ i.e.,
$$(\sigma_{\mu},\sigma_{\hat{\nu}})_q=\delta_{\mu,\nu}.$$
\end{corollary}

\section{TOTAL POSITIVITY}\label{sec:8}
A matrix $A\in GL_{n}(\R)$ is said to be $totally$ $positive$ $(
\textrm{resp}.\hs totally\hs nonnegative)$ if all the minors of
$A$ are positive (resp. nonnegative). These matrices form
semialgebraic subset of $GL_n(\R).$ Total positivity for
$GL_n(\R)$ was mainly studied around the 1950's by Schoenberg,
Gantmacher-Krein and others, and has diverse applications such as
oscillating mechanical systems and planar Markov process. Since it
was generalized to all the reductive algebraic groups by Lusztig
in the early 1990's \cite{Lus1}, total positivity has been more
noted by its connections with the canonical bases. Furthermore,
very recently it was shown by Rietsch that in type $A,$  total
positivity of $\mathcal{V}_P$  has a very close connection with
the positivity of Schubert basis functions on $\mathcal{V}_P$,
more precisely the set of the totally nonnegative elements in the
varieties $\mathcal{V}_P$ coincides with the set where  the
Schubert basis functions are nonnegative (\cite{Riet1}), and it
was conjectured that this will be true for the other classical Lie
types (\cite{Riet2}). In this section, we explicitly describe the
totally positive parts of the varieties $\mcvc$ and $\mcvb,$ and
characterize them via the Schubert basis functions.
\subsection{ Total positivity of $U^+$ and its Bruhat
cells}\label{total posi} Let $U^+(\R)$ be the group of real points
in $U^+.$ For each $i=1,...,n,$  let $x_i(t):=\textrm{exp}(te_i)$
be the one parameter subgroup corresponding to the simple root
vector $e_i.$ Define the subset $U^{+}(\R_{ \geq 0})$ of $totally$
$nonnegative$ $elements$ in $U^+$ as the multiplicative semigroup
generated by all $x_i(t)$ with $i=1,...,n,$ and $t\geq 0.$ For any
$w\in W,$ we set $U_w^{+}(\R _{>0}):=U^{+}(\R_{\geq 0})\cap
B^-wB^-.$ Then since $G=\bigcup_{w\in W}B^-wB^-,$ we have
$U^{+}(\R _{\geq 0})=\bigcup_{w\in W}U_w^{+}(\R_{>0}).$  The
algebraic  subset $U_{w_0}^{+}(\R_{>0})$ is called the $totally$
$positive$ $part$ of $U^+$, denoted by $U^{+}(\R_{>0}),$
 and  it is open in $U^{+}(\R)$ and its closure is the set of all totally nonnegative elements
 in $U^+, $ i.e., $\overline{U^{+}(\R_{>0})}=U^+(\R_{\geq 0}).$ If   $ G$ is $GL_n(\R),$
 the above definitions coincide  with the classical ones.
In case $G$ is $Sp_{2n}(\C)$ (resp. $SO_{2n+1}(\C)$), the total
positivity of G can be understood  via that of
$\widetilde{G}:=SL_{2n}(\C)$ (resp. $SL_{2n+1}(\C))$. With the
matrices $J$ in $\ref{basic:sp}$ and \ref{subsec 2n+1} for
$Sp_{2n}(\C)$ and $SO_{2n+1}(\C)$, respectively,  the groups $G$
is naturally embedded into $\widetilde{G}.$  More precisely, we
have \be \label{total preserve}U^+=\widetilde{U}^+ \cap G, \hs
\textrm{and}\hs T=\widetilde{T}\cap G,\ee where $\widetilde{U}^+$
and $\widetilde{T}$ are unipotent upper triangular matrices and
diagonal matrices of $SL_{2n}(\C)$ or $SL_{2n+1}(\C)$.
Furthermore, this embedding preserves the total positivity of
$U^+$ and Bruhat cells of $U^+$(\cite{BZ1}), that is, for each
$w\in W,$ \be\label{total positivity via type A} U_{w}^{+}(\R
_{>0})=\widetilde{U}_{\tilde{w}}^{+}(\R _{>0})\cap G, \hs
\textrm{and}\hs U^{+}(\R _{\geq 0})=\widetilde{U}^{+}(\R _{\geq
0})\cap G,\ee where $\tilde{w}$ is the image of $w$ under a proper
embedding (\cite{BZ1})of $W$ into $\widetilde{W}=S_{2n}.$
 \begin{definition}
For each parabolic subgroup $P \subset G,$ denote by
$\mathcal{V}_P(\R)$ the set of all the real points of
$\mathcal{V}_P \subset G^\vee.$ Define
$$\mathcal{V}_P(\R_{\geq 0}):=\mathcal{V}_P(\R)\cap (U^{\vee})^+(\R_{\geq0}),
\hs \textrm{and}\hs \mathcal{V}_P(\R_{>0}):=\mathcal{V}_P(\R_{\geq
0})\cap (B^\vee)^-w^P (B^\vee)^-.$$
\end{definition}

 Note that  if $P$ is a maximal parabolic
subgroup, then \be \label{break-up}
\mathcal{V}_P=((U^\vee)^+)^{e^\vee}\cap\overline{ (B^\vee)^-w^P
(B^\vee)^-}=[((U^\vee)^+)^{e^\vee}\cap (B^\vee)^-w^P(B^\vee)^-]
\cup \{id\},\ee and hence we have \be \label{break-up3}
\mathcal{V}_P(\R_{\geq 0})=\mathcal{V}_P(\R_{>0}) \cup \{id\}.\ee

 When $P$ is the maximal parabolic subgroup  $P_n$ of $SO_{2n+1}(\C)$ and  $Sp_{2n}(\C),$ or  $P_{n+1}$ of
 $SO_{2n+2}(\C),$
 to be consistent with the previous notations,  we write, for example, $\mcvc(\R_{\geq 0})$
 rather than
 $\mathcal{V}_{P_n}(\R_{\geq0})$ for the triple $(SO_{2n+1}(\C), P_n, Sp_{2n}(\C)).$
\subsection{Total positivity in Lie type A}\label{subsec:Total A}
By (\ref{total positivity via type A}), we can detect totally
nonnegative elements in the varieties $\mcvc$ and $\mcvb$ by
embedding these varieties to a `relevant' variety in Lie type $A$,
whose total positivity is well understood. Now we define the
`relevant' variety in Lie type $A$ and describe its total
positivity. Consider the set $\mathcal{V}_n:=\{u_{2n}(t\zeta^I)\in
SL_{2n}(\C)\mid t \in \C, I \in \mathcal{T}_n\}$, with
$u_{2n}(t\zeta^I)$ defined in $\ref{def and nota}.$ This set can
be realized as a variety $\mathcal{V}_{P_n}$ whose coordinate ring
is isomorphic to the quantum cohomology ring of Grassmannian
manifold $SL_{2n}(\C)/P_n,$ where $P_n$ is the maximal parabolic
subgroup corresponding to the $n$-th fundamental weight. Then
Theorem $8.4$ in \cite{Riet1} implies that the set of totally
nonnegative elements in $\mathcal{V}_n$ is
$\mathcal{V}_n(\R_{\geq0}):=\{u_{2n}(t\zeta^{I_0}) \mid t\geq 0
\}.$ Letting $\mathcal{V}_n(\R_{>0}):=\{u_{2n}(t\zeta^{I_0})\mid
t\in \R_{>0}\}$, we have  \be \label{break up 2}
\mathcal{V}_n(\R_{\geq 0})=\mathcal{V}_n(\R_{>0})\cup \{id\}.\ee
Therefore we get
$\mathcal{V}_n(\R_{\geq0})=\mathcal{V}_{P_n}(\R_{\geq0})$ by
(\ref{break-up3}).

On the other hand, Proposition $9.3$ in \cite{Riet1}  implies that
the set $\mathcal{V}_n(\R_{>0})$ is characterized by the
positivity of `special' Schur polynomials, more precisely, \be
\label{criterion} \mathcal{V}_n(\R_{>0})=\{u_{2n}(t\zeta^I)\in
\mathcal{V}_{n}\mid S_{(m^k)}(t\zeta^I)>0\hspace{0.05in}\textrm{
for all} \hspace{0.05in} m,k\leq n\}.\ee

\subsection{Total positivity in $\mcvc$ and $\mcvb$}
From \ref{total posi} we now have `natural' embeddings  $\mcvc
\hookrightarrow \mathcal{V}_n \hookrightarrow SL_{2n}(\C)$ and
$\mcvb \hookrightarrow SL_{2n+1}(\C),$ and so we will consider
elements of $\mcvc$ (resp. $\mcvb$) as elements of $SL_{2n}(\C)$
(resp. $SL_{2n+1}(\C)$). Note that while $\mcvc$ has an
`intermediate' variety $\mathcal{V}_n$, $\mcvb$ does not. But an
element $u_{2n+1}(t\zeta^I)$ of $\mcvb \subset SL_{2n+1}$ can be
viewed as a submatrix of $u_{2n+2}(t\zeta^I)\in \mathcal{V}_{n+1}
\subset SL_{2n+2}(\C),$ by taking the first $(2n+1)$-rows and
columns  of $u_{2n+2}(t\zeta^I)$.  So we can use the total
positivity of $\mathcal{V}_{n+1}$ to describe that of $\mcvb.$
\begin{theorem} \label{thm:total positivity 1}  The totally
nonnegative elements in $\mcvc$ can be described as follow. \bn
\item $\mcvc(\R_{\geq 0})=\{u_{2n}(t\zeta^{I_0})\mid t \in \R_{\geq 0}\}=
\{u\in \mcvc \mid \dot{\Pt}_{\lambda}(u)\geq0
\hspace{0.05in}\textrm{for all} \hspace{0.05in}\lambda \in \dn
\},$
\item $\mcvc(\R_{>0})=\{u_{2n}(t\zeta^{I_0})\mid t \in \R_{>0}\}=
\{u\in \mcvc\mid \dot{\Pt}_{\lambda}(u)>0
\hspace{0.05in}\textrm{for all} \hspace{0.05in}\lambda \in \dn
\}.$ \en
\end{theorem}

\begin{proof}
Consider the embedding $\mcvc \hookrightarrow \mathcal{V}_n $.
Then the first identity of $(1)$ follows from the result in type
$A$ of $\ref{subsec:Total A}$ and $(\ref{total positivity via type
A})$.  The second identity of $(1)$ is a direct consequence of
$(1)$ of  Theorem \ref{positivity of Q}.  We use the first
identity of $(1)$ and (\ref{break-up3}) to get the first identity
of $(2)$. The second identity of $(2)$ follows from $(1)$ of
Theorem \ref{positivity of Q} and the fact that for any $I\in
\mathcal{I}_n $ $\Pt_\lambda(t\zeta^I)=0$ for all $\lambda\in
\mathcal{D}(n)$ if and only if $t=0$.
\end{proof}
The similar analysis can be applied to $G=SO_{2n+1}(\C)$ and
$\widetilde{G}=SL_{2n+1}(\C).$

\begin{theorem}\label{to2}
The totally nonnegative elements in $\mcvb$ can be described as
follow. \bn
\item $\mcvb(\R_{\geq 0})=\{u_{2n+1}(t\zeta^{I_0})\mid t \in \R_{\geq 0}\}=
\{u\in \mcvb \mid \dot{\Qt}_{\lambda}(u)\geq0
\hspace{0.05in}\textrm{for all} \hspace{0.05in}\lambda \in \dn
\},$
\item $\mcvb(\R_{>0})=\{u_{2n+1}(t\zeta^{I_0}) \mid t \in \R_{>0}\}=
\{u\in \mcvb \mid \dot{\Qt}_{\lambda}(u)>0
\hspace{0.05in}\textrm{for all} \hspace{0.05in}\lambda \in \dn
\}.$ \en
\end{theorem}
\begin{proof}
Since $\mathcal{V}_{n+1}(\R_{\geq 0})=\{
u_{2n+2}(t\zeta^{I_0})\mid t\in \R_{\geq 0}\},$ and since if
$u_{2n+2}(t\zeta^I)$ is totally nonnegative, so is the submatrix
$u_{2n+1}(t \zeta^I)$, $u_{2n+1}(t\zeta^{I_0})$ for $t \geq0$ are
totally nonnegative elements of $\mcvb,$ i.e.,
$$\{u_{2n+1}(t \zeta^{I_0})\mid t \in \R_{\geq 0}\} \subseteq
\mcvb(\R_{\geq 0}).$$ For the converse inclusion, considering
$(\ref{break-up3}),$ it suffices to show that
$$\mcvb(\R_{>0})\subseteq \{u_{2n+1}(t \zeta^{I_0})\mid t \in \R_{>
0}\}.$$ So let $u \in \mcvb(\R_{>0})$. We can write
$u=u_{2n+1}(t_{0}\zeta^{I})$ for some nonzero $t_{0}\in \C$ and $I
\in \mathcal{I}_{n+1}.$ Take $\bar{u}=u_{2n+2}(t_0\zeta^I) \in
\mathcal{V}_{n+1}.$ We claim that $S_{(m^k)}(t_0\zeta^I)>0$ for
all $m, n \leq n+1.$ For $m,n $ with $m\ne n+1$, or $k\ne n+1,$
$S_{(m^k)}(t_0\zeta^I)$ is the determinant of a submatrix $v$ of
$u_{2n+1}(t_{0}\zeta^I)$, which is nonnegative since
$u_{2n+1}(t_0\zeta^I)$ is totally nonnegative, but this cannot be
zero since the row vectors of $v$ are linearly independent. For
$m,k=n+1,$ $S_{(m^k)}(t_0\zeta^I)=(a_{n+1})^{n+1}$, which is
positive since $a_{n+1}=E_{n+1}(t_0 \zeta^I)$ is an nonzero entry
of totally nonnegative element $u_{2n+1}(t_0\zeta^I)$. Therefore
this satisfies the criterion $(\ref{criterion})$ for $\bar{u} \in
\mathcal{V}_{n+1}$ to be totally nonnegative elements. But the
totally nonnegative elements of $\mathcal{V}_{n+1}$ are
$u_{2n+2}(t\zeta^{I_0})$ with $t\in \R_{\geq0},$ so we have
$I=I_0$ and $t_0 >0.$ Therefore we have
$$\mcvb(\R_{> 0}) \subseteq \{u_{2n+1}(t \zeta^{I_0})\mid t \in
\R_{> 0}\}.$$ This proves the first identity of $(1)$. For the
second identity of $(1),$ as in Theorem \ref{thm:total positivity
1}, we use Theorem \ref{positivity of Q} to  get $$\{u_{2n+1}(t
\zeta^{I_0})\mid t \in \R_{\geq 0}\}=\{u\in \mcvb \mid
\dot{\Qt}_{\lambda}(u)\geq0 \hspace{0.05in}\textrm{for all}
\hspace{0.05in}\lambda \in \mathcal{D}(n+1)\},$$ which is a subset
of $$\{u\in \mcvb \mid \dot{\Qt}_{\lambda}(u)\geq0
\hspace{0.05in}\textrm{for all} \hspace{0.05in}\lambda \in
\mathcal{D}(n)\}.$$  So it is enough to show that $$\{u\in \mcvb
\mid \dot{\Qt}_{\lambda}(u)\geq0 \hspace{0.05in}\textrm{for all}
\hspace{0.05in}\lambda \in \mathcal{D}(n)\}\subseteq \{u\in \mcvb
\mid \dot{\Qt}_{\lambda}(u)\geq0 \hspace{0.05in}\textrm{for all}
\hspace{0.05in}\lambda \in \mathcal{D}(n+1)\}.$$ So suppose $u$ is
an element of $\mcvb$ such that $\dot{\Qt}_\lambda(u)\geq0$ for
all $\lambda \in \dn.$ Note that  the relation
$X_{n}^2=2X_{n-1}X_{n+1}$ implies  the relation,
$$(\dot{\Qt}_{n}(u))^2=2 \dot{\Qt}_{n-1}(u)\dot{\Qt}_{n+1}(u).$$
Since $\dot{\Qt}_{n}(u)$ and $\dot{\Qt}_{n-1}(u)$ are nonnegative,
so is $\dot{\Qt}_{n+1}(u).$ Let $\lambda \in \mathcal{D}(n+1).$ If
$\lambda\in \dn,$ by the assumption, $\dot{\Qt}_\lambda(u)$ is
nonnegative. If $\lambda \notin \dn,$ we can write $\lambda=(n+1)
\cup \mu$ for some $\mu\in \dn.$ If we apply the factorization of
$\Qt$-polynomials, we have
$$\dot{\Qt}_\lambda(u)=\dot{\Qt}_{n+1}(u) \dot{\Qt}_\mu(u).$$
Since $\dot{\Qt}_{n+1}(u)$ and $ \dot{\Qt}_\mu(u)$ are
nonnegative, so is $\dot{\Qt}_\lambda(u),$ as desired. The
identities of $(2)$ can be obtained by doing the same argument as
in the proof of $(2)$ in Theorem \ref{thm:total positivity 1}.
\end{proof}

From Theorem $\ref{thm:total positivity 1}$ and Theorem
$\ref{to2},$ we conclude that Rietsch's conjecture about the
relation between the total positivity and quantum cohomology ring
holds for the triples $(SO_{2n+1}(\C),Sp_{2n}(\C), P_n)$ and
$(Sp_{2n}(\C),SO_{2n+1}(\C),P_n).$

{\it Acknowledgements.}  The author would like to thank K. Rietsch
for useful correspondence about Peterson's results, A. Zelevinsky
for answering many questions about the total posivity over
e-mails, S. Seo for useful discussion about combinatorial aspects,
and H. Tamvakis for valuable suggestions and discussions.


\begin{thebibliography}{50}
\bibitem{BKT} A. Buch  A. Kresch and H. Tamvakis, {\it Gromov-Witten
invariants on Grassmannians}, J. Amer. Math. Soc. 16(2003),
901-915.
\bibitem{BZ1} A. Berenstein and A. Zelevinsky, {\it Total Positivity in Schubert Varieties}, Comment. Math. Helv.
72(1997), 128-166.
\bibitem{FH}W. Fulton and J. Harris, {\em Representation Theory}, Springer-Verlag, New York 1991.
\bibitem{HB} H. Hiller and B. Boe, {\it Pieri formula for $SO_{2n+1}/U_n$
and  $Sp_n/U_n$}, Adv. in Math. $62$ (1986), 49-67.
\bibitem{Ko} B. Kostant, {\it Flag manifolds quantum cohomology, the
Toda lattice, and the representation with highest weight $\rho$},
Selecta Math.(N.S.) $2$ $(1996)$, 43-91.
\bibitem{Knu}D. Knutson, {\em $\lambda$-rings and the representation
theory of the symmetric group}, Springer Lecture Notes 308, 1973.
\bibitem{KT2}  A. Kresch and H. Tamvakis, {\it Quantum cohomology of the
Lagrangian Grassmannian}, J. Algebraic Geometry 12 (2003),
777-810.
\bibitem{KT1} A. Kresch and H. Tamvakis, {\it Quantum cohomology of
orthogonal Grassmannians}, Compositio Math. 140 (2004), 482-500.
\bibitem{Las and Pra} A. Lascoux and P. Pragacz, {\it Operator
Calculus for $\Qt$-Polynomials and Schubert Polynomials}, Adv. in
Math $140$ (1998), 1-43.
\bibitem{Lus1} G. Lusztig, {\em Total positivity in reductive groups, Lie
theory and geometry} : in honor of Bertram Kostant(G. I.
Lehrer,ed.), Progress in Mathematics, vol. 123, Birkhauser,
Boston,1994, pp. 531-568.
\bibitem{Lus2} G. Lusztig, {\em Total positivity and canonical bases,
Algebraic groups and Lie groups}, Cambridge Univ. Press,
Cambridge, 1997, pp.281-295.
\bibitem{Mac1} I.G. Macdonald, {\em Symmetric Functions and Hall Polynomials}, Second
edition, Oxford Univ. Press, 1995.
\bibitem{Minc} H. Minc, {\em Nonnegative Matrices}, John Willy \&
Sons, 1988.
\bibitem{Pete1} D. Peterson, {\em Qunatum Cohomology of $G/P$}, Lecture
Course, spring term, M.I.T., 1997.
\bibitem{Pete2} D. Peterson, {\em Qunatum Cohomology of $G/P$}, Seminaire
de Mathamatiques Superieures: Representation Theories and
Algebraic Geometry, Universite de Montreal, Canada, July 28-Aug.
8, 1997(unpublished lecture notes).
\bibitem{Pra and Rat 1} P. Pragacz and J. Ratajski, {\it Formulas for Lagrangian and orthogonal
degeneracy loci; $\Qt$-polynomial approach}, Compositio
Mathematica 107; 11-87, 1997.
\bibitem{Riet1} K. Rietsch,  {\it Quantum cohomology rings of
Grassmannians and total positivity}, Duke Mathematical Journal,
Vol 110, no. 3 (2001), 523-553.
\bibitem{Riet2} K. Rietsch, {\it Totally positive Toeplitz matrices and quantum cohomology of
partial flag varieties}, J. Amer. Math. Soc, Vol. 16(2003),
363-392
\bibitem{Riet3} K. Rietsch, {\it A mirror symmetric construction of
$H_T^*(G/P)_{(q)}$}, arXiv: math.AG/0511124.



\end{thebibliography}
\end{document}